\documentclass[11pt]{article}
\usepackage{mathrsfs}
\usepackage[centertags]{amsmath}
\usepackage{amsfonts}
\usepackage{amssymb}
\usepackage{amsthm}
\usepackage{indentfirst}
\usepackage{epsfig}
\usepackage {graphicx}
\usepackage{color}

\date{}

\definecolor{c20}{rgb}{0.,0.7,0.}
\definecolor{c30}{rgb}{0.98,0.00,0.00}
\definecolor{c40}{rgb}{1,0.1,0.7}
\definecolor{c50}{rgb}{1,0,0}

\def\pzx#1{\textcolor{c50}{#1}}
\def\pzx#1{#1}
\newtheorem{theorem}{Theorem}[section]

\newtheorem{lemma}{Lemma}[section]

\newtheorem{remark}{Remark}[section]

\numberwithin{equation}{section}

\def\P{\operatorname*{\mathbf{P}}}
\def\E{\operatorname*{\mathbf{E}}}

\def\Cov{\operatorname*{\mathbf{Cov}}}

\def\M{\operatorname*{\mathbf{M}}}
\def\O{\operatorname*{\mathbf{O}}}
\def\T{\operatorname*{\mathbf{T}}}

\def\q{\operatorname*{\mathbf{q}}}
\def\t{\operatorname*{\mathbf{t}}}
\def\a{\operatorname*{\mathbf{a}}}
\def\s{\operatorname*{\mathbf{s}}}
\def\S{\operatorname*{\mathbf{S}}}
\def\k{\operatorname*{\mathbf{k}}}
\def\u{\operatorname*{\mathbf{u}}}
\def\v{\operatorname*{\mathbf{v}}}
\def\l{\operatorname*{\mathbf{l}}}
\def\I{\operatorname*{\mathbf{I}}}
\def\i{\operatorname*{\mathbf{i}}}
\def\j{\operatorname*{\mathbf{j}}}
\def\m{\operatorname*{\mathbf{m}}}
\allowdisplaybreaks[4]

\begin{document}
\title{Maxima and minima of homogeneous Gaussian random fields over continuous time and uniform grids}
\author{{Yingyin Lu \;\;\; Zuoxiang Peng\thanks{Corresponding author. Email: pzx@swu.edu.cn} }    \\
{\small School of Mathematics and Statistics, Southwest University, Chongqing, 400715, China}}

\maketitle
\begin{quote}
{\bf Abstract.}~~In this paper, for centered homogeneous Gaussian random fields  the joint limiting distributions of normalized maxima and minima over continuous time and uniform grids are investigated. It is shown that maxima and minima are asymptotic dependent for strongly dependent homogeneous Gaussian random field with the choice of sparse grid, Pickands' grid or dense grid,  while for the weakly dependent Gaussian random field maxima and minima are asymptotically independent.


\end{quote}

\section{Introduction}\label{sec1}
\label{sec1}
Let $\{ X(t), t\geq 0\}$ be a stationary Gaussian process with mean zero, variance one and the correlation function $r(t)$ \pzx{which satisfies for} some $\alpha \in (0,2]$,
\begin{eqnarray}
\label{eq1.1}
r(t)=1-|t|^{\alpha}+o\left(|t|^{\alpha}\right),\;\ t\to 0.
\end{eqnarray}
Assume that
\begin{eqnarray}
\label{eq1.2}
r(t)\log t\to 0\;\ \mbox{as}\;\ t\to \infty.
\end{eqnarray}
Under conditions \eqref{eq1.1} and \eqref{eq1.2}, the limit distribution theory on the maximum of $\{ X(t), t\geq 0\}$ up to time $T$,
\[ M_{T}=\max\{ X(t) , t\in [0,T]\}\]
is well developed that
\begin{eqnarray*}
\P\{ a_{T}\left(M_{T}-b_{T}\right)\leq x\} \to \Lambda(x)
\end{eqnarray*}
as $T\to\infty$, where $ \Lambda(x)= \exp\left(-e^{-x}\right) $ is the standard \pzx{Gumbel} distribution and
\begin{eqnarray*}
a_{T}=\sqrt{2 \log T},\;\ b_{T}=\sqrt{2 \log T} + \frac{\log\Big[(2\pi)^{-1/2} H_{\alpha} (2\log T)^{-1/2+1/\alpha}\Big]}{\sqrt{2 \log T}}.
\end{eqnarray*}
Here $ H_{\alpha} $ is the well-known Pickands' constant, which is defined by $ H_{\alpha}=\lim_{\lambda\to\infty} H_{\alpha}(\lambda)/\lambda $ with
\begin{eqnarray}
\label{eq.1.3}
H_{\alpha}(\lambda)=\E\exp\left( \max_{t\in [0,\lambda]} \sqrt{2} B_{\alpha/2}(t)-t^{\alpha}\right)
\end{eqnarray}
and $ B_{H}$ is a fractional Brownian motion with $ \E B_{H}^{2}(t)=|t|^{2H}$. It is well known that $ 0<H_{\alpha}<\infty $, \pzx{see e.g.} Leadbetter et al. (1983), Pickands (1969) and Piterbarg (1996). The limit distribution theorem about $ M_{T}$ extended by Mittal and Ylvisaker (1975) and McCormick and Qi (2000), and D\c{e}bicki et al. (2013) extended the results to homogeneous Gaussian random fields.

\pzx{It is known that extreme value theory of Gaussian random fields  may be applied to image analysis, quantum chaos, queuing theory, insurance mathematics, number theory and so on; see e.g. Adler (2000), Adler et al. (2014), and Hashorva and Ji (2016). Further, numerical simulation of trajectories of high extremes of continuous random processes may be performed through the discrete time random processes depending heavily on the sampling frequency, see, for instance Leadbetter et al. (1983), Piterbarg (2004) and recent work of Song et al. (2018, 2019) on flaw detection by using ultrasonic response signals.}

The joint limiting distributions of $M_{T}$ and its discrete time maximum $M_{T}^{\delta}=\max\{ X(t),t\in [0,T]\cap \Re(\delta)\}$ was first studied by Piterbarg (2004) under the conditions \eqref{eq1.1} and \eqref{eq1.2}, where the uniform \pzx{grid} is given by $\Re=\Re(\delta)=\{k\delta,\;\ k=0,1,2\ldots\}$ with $\delta \left( 2\log T \right)^{1/\alpha} \to D\in[0,\infty] $  as $T\to \infty.$  For more details, see Piterbarg (2004). The results \pzx{on} multivariate stationary Gaussian processes can be found in Tan and Hashorva (2014, 2015), Tan and Wang (2013) and Tan and Tang (2014) for strongly dependent Gaussian processes. Further,  Turkman (2012) considered the non-Gaussian processes and Chen and Tan (2016) studied the asymptotic \pzx{behavior} of $M_{T}$, $M_{T}^{\delta}$ and the partial sum of dependent Gaussian processes.

For the joint asymptotic behaviors of $M_{T}$ and $M_{T}^{\delta}$ of homogeneous Gaussian random field,  Tan and Wang (2015) considered the following model. Let $\{ X(\t), \t\geq \mathbf{0}\}$ be a homogeneous Gaussian field with mean zero, variance one and covariance function $ r(\mathbf{t})=\Cov(X(\mathbf{t}),X(\mathbf{0})) $ satisfies the following conditions, for $ d\geq 2 $:

$\mathbf{A1}: r(\t)=1-\sum_{i=1}^{d}|t_{i}|^{\alpha_{i}}(1+o(1))$, as $t_{i}\to 0$, with $\alpha_{i}\in (0,2] $ ;

$\mathbf{A2}: r(\t) < 1$, for $\t \neq \mathbf{0} $ ;

$\mathbf{A3}: \lim_{\T\to \infty} r(\T)\log(\prod_{i=1}^{d} T_{i})=r\in[0,\infty)$, where $\T\to \infty$ means $T_{i}\to\infty, i=1,2,\cdots,d$. If $T_{i}=0$, $r(T_{1},\cdots,T_{i-1},0,T_{i+1},\cdots, T_{d})\log(\prod_{j\neq i}^{d}T_{j})$ is bounded.

Define
\[ \mathbf{M}_{\T}=\left( M_{\T}, M_{\T}^{\boldsymbol{\delta}}\right)=\left( \max_{\t\in \I_{\T}} X(\t),\max_{\t\in \I_{\T}\bigcap \prod_{i=1}^{d} \Re(\delta_{i})} X(\t) \right),\]
where $ \I_{\T}=\prod_{i=1}^{d}[0,T_{i}]$, and $ \I_{\T}\bigcap \prod_{i=1}^{d} \Re(\delta_{i}) $ means  $ \prod_{i=1}^{d} \{[0,T_{i}] \bigcap \Re(\delta_{i})\} $, and the uniform \pzx{grid}
$ \Re(\delta_{i})= \{k \delta_{i}, k\in N\} $ is given by
\[ \delta_{i}\left( 2\log\prod_{i=1}^{d} T_{i}\right)^{1/\alpha_{i}} \to D_{i},\;\ i=1,2,\cdots,d.\]
\pzx{We say that the grid is dense if all $ D_{i}=0 $} and if all $ D_{i}=\infty $, the grid is sparse. The grid is called a Pickands grid if all $ D_{i}\in (0,\infty) $. Under conditions $\mathbf{A1}-\mathbf{A3}$, Tan and Wang (2015) derived the limiting distribution of $ \M_{\T} $ when the uniform grid is sparse grid, Pickands' grid and dense grid, respectively.

 In this paper, our focus is on the joint limit distributions of maxima and minima of aforementioned homogeneous Gaussian random fields.  Davis (1979) established the joint limiting distribution of maxima and minima of weakly dependent stationary sequences, and weakly dependent stationary Gaussian processes was studied by Berman (1971). For the asymptotic distributions of maxima and minima on bivariate H\"{u}sler-Reiss models, see Liao and Peng (2015) and Lu and Peng (2017).

\pzx{Similarly} to the definition of maxima $ \M_{\T} $, define the minima $ \m_{\T}$ as follows:
\[ \mathbf{m}_{\T}=\left( m_{\T}, m_{\T}^{\boldsymbol{\delta}}\right)=\left( \min_{\t\in \I_{\T}} X(\t),\min_{\t\in \I_{\T}\bigcap \prod_{i=1}^{d} \Re(\delta_{i})} X(\t) \right)\]
and  let
\begin{eqnarray*}
\begin{cases}
b_{\T}^{\boldsymbol{\delta}}=a_{\T}+a_{\T}^{-1}\log\left((2\pi)^{-1/2}(\prod_{i=1}^{d} \delta_{i}^{-1})(a_{\T})^{-1}\right)\\
b_{\a,\T}=a_{\T}+a_{\T}^{-1}\log\left((2\pi)^{-1/2}(\prod_{i=1}^{d} H_{a_{i},\alpha_{i}})(a_{\T})^{\sum_{i=1}^{d} \frac{2}{\alpha_{i}}-1}\right)\\
b_{\T}=a_{\T}+a_{\T}^{-1}\log\left((2\pi)^{-1/2}(\prod_{i=1}^{d} H_{\alpha_{i}})(a_{\T})^{\sum_{i=1}^{d} \frac{2}{\alpha_{i}}-1}\right),
\end{cases}
\end{eqnarray*}
where $a_{\T}=\sqrt{2\log(\prod_{i=1}^{d} T_{i})}$ and $H_{a_{i},\alpha_{i}}=\lim_{\lambda_{i}\to\infty}\frac{H_{a_{i},\alpha_{i}}(\lambda_{i})}{\lambda_{i}}\in(0,\infty)$
with
\[H_{a_{i},\alpha_{i}}(\lambda_{i})=\E\exp\left( \max_{k_{i}a_{i}\in[0,\lambda_{i}]} \sqrt{2} B_{\alpha_{i}/2}^{(i)}(k_{i}a_{i})-(k_{i}a_{i})^{\alpha_{i}}\right),\; i=1,2,\cdots,d,\]
 where $B_{\alpha_{1}/2}^{(1)}(\cdot),\cdots, B_{\alpha_{d}/2}^{(d)}(\cdot)$ are independent fractional Brownian motions, cf. Piterbarg (2004) and Tan and Wang (2015). Further, let the bivariate normalizing constants $\mathbf{u}_{\T}$ and $\mathbf{v}_{\T}$ be given by
\begin{eqnarray}\label{eq4.2}
\begin{cases}
 \mathbf{u}_{\T}=\left( u_{\T}(x_{2}), u_{\T}^{\boldsymbol{\delta}}(y_{2})\right)=\left(b_{\T}+\frac{x_{2}}{a_{\T}},b_{\T}^{*}+\frac{y_{2}}{a_{\T}}\right) \\
 \mathbf{v}_{\T}=\left( v_{\T}(x_{1}), v_{\T}^{\boldsymbol{\delta}}(y_{1})\right)=\left(-b_{\T}+\frac{x_{1}}{a_{\T}},-b_{\T}^{*}+\frac{y_{1}}{a_{\T}}\right),
 \end{cases}
 \end{eqnarray}
where $ b_{\T}^{*}=b_{\T}^{\boldsymbol{\delta}}$ for sparse grids, $ b_{\T}^{*}=b_{\a,\T} $ for Pickands grids and $ b_{\T}^{*}=b_{\T} $ for dense grids.

Throughout this paper, let $ \phi(x) $ and $ \Phi(x) $ denote respectively the density function and distribution function of a standard normal random variable, and $ \Psi(x)=1-\Phi(x) $, and operations of vectors \pzx{mean} componentwise operating. For example, for vectors $ \t=(t_{1},t_{2},\cdots,t_{d}) $ and $ \s=(s_{1},s_{2},\cdots,s_{d}) $, operations of $ \s\leq \t $, $\s-\t$, $ \s \t$, and $ \s ^{\t}$ mean $ s_{i}\leq t_{i}, i=1,2,\cdots,d $, $(s_{1}-t_{1},s_{2}-t_{2},\cdots,s_{d}-t_{d})$, $ (s_{1}t_{1},s_{2}t_{2},\cdots,s_{d}t_{d}) $ , and $ (s_{1}^{t_{1}},s_{2}^{t_{2}},\cdots,s_{d}^{t_{d}}) $, respectively. Let $C$ be positive constant with values varying from place to place.

The contents of this paper are organized as follows. Section \ref{sec2} presents the main results and Section \ref{sec3} gives some auxiliary lemmas. The proofs of the main results will be given in Section \ref{sec4}.

\section{Main results}\label{sec2}

\begin{theorem}
\label{th1}
Let $X(\t)$ be a centered homogenous Gaussian field with \pzx{unit variance} and covariance function $r(\t)$ satisfying $\mathbf{A1}-\mathbf{A3}$. Then for any sparse grids $ \Re(\delta_{i}) $, $ i=1,2,\cdots,d $,
\begin{eqnarray*}
&&\lim_{\T\to\infty}\P\left(\mathbf M_{\T}\leq\mathbf u_{\T},\mathbf m_{\T}\leq\mathbf v_{\T}\right)\nonumber\\
&=&\int_{-\infty}^{+\infty}\Big[1-\exp\left(-e^{x_{1}+r+\sqrt{2r}z}\right)-\exp\left(-e^{y_{1}+r+\sqrt{2r}z}\right)
+\exp\left(-e^{x_{1}+r+\sqrt{2r}z}-e^{y_{1}+r+\sqrt{2r}z}\right)\Big]\nonumber\\
&&\times\exp\Big\{-\left(e^{-x_{2}-r+\sqrt{2r}z}+e^{-y_{2}-r+\sqrt{2r}z}\right)\Big\}d\Phi(z),
\end{eqnarray*}
where $\mathbf{u}_{\T}$, $\mathbf{v}_{\T}$ are given by \eqref{eq4.2}.
\end{theorem}

\pzx{Define,
\begin{eqnarray*}
H_{\a,\boldsymbol{\alpha}}^{x,y}(\boldsymbol{\lambda})=\int_{-\infty}^{\infty} e^{s}\P\Big( \max_{\k\a\in\prod_{i=1}^{d}[0,\lambda_{i}]}\sqrt{2}\chi(\k\a)>s+x,  \max_{\t\in\prod_{i=1}^{d}[0,\lambda_{i}]} \sqrt{2}\chi(\t)>s+y \Big) ds,
\end{eqnarray*}
where $\chi(\t)=\sum_{i=1}^{d} B_{\alpha_{i}/2}^{(i)}(t_{i})- \sum_{i=1}^{d}|t_{i}|^{\alpha_{i}}$.}

\begin{theorem}
\label{th2}
Let $X(\t)$ be a centered homogenous Gaussian field with \pzx{unit variance} and covariance function $r(\t)$ satisfying $\mathbf{A1}-\mathbf{A3}$ and $\mathbf{u}_{\T}$, $\mathbf{v}_{\T}$ be given by \eqref{eq4.2}. \pzx{Then for any Pickands grids $\Re(a_{i}(2\log \prod_{i=1}^{d} T_{i})^{-1/\alpha_{i}})$ with $a_{i}>0$,  $ i=1,2,\cdots,d $, the following limit exists,
\[ H_{\a,\boldsymbol{\alpha}}^{x,y}=\lim_{\boldsymbol{\lambda}\to\infty} H_{\a,\boldsymbol{\alpha}}^{x,y}(\boldsymbol{\lambda})/(\prod_{i=1}^{d} \lambda_{i}) \in (0,\infty)  \]}
and
\begin{eqnarray*}
&&\lim_{\T\to\infty}\P\left(\mathbf M_{\T}\leq\mathbf u_{\T},\mathbf m_{\T}\leq\mathbf v_{\T}\right)\\
&=&\int_{-\infty}^{+\infty}\Big[1-\exp\Big(-e^{x_{1}+r+\sqrt{2r}z}\Big)-\exp\Big(-e^{y_{1}+r+\sqrt{2r}z}\Big)
+ \exp\Big(-e^{x_{1}+r+\sqrt{2r}z}-e^{y_{1}+r+\sqrt{2r}z}\\
&&-H_{\a,\boldsymbol{\alpha}}^{-x_{1}+\log( \prod_{i=1}^{d}H_{\alpha_{i}}),-y_{1}+\log(\prod_{i=1}^{d} H_{a_{i},\alpha_{i}})}e^{r+\sqrt{2r}z}\Big)\Big]
\times \exp\Big\{-\Big(e^{-x_{2}-r+\sqrt{2r}z}+e^{-y_{2}-r+\sqrt{2r}z}\\
&&-H_{\a,\boldsymbol{\alpha}}^{x_{2}+\log( \prod_{i=1}^{d}H_{\alpha_{i}}),y_{2}+\log(\prod_{i=1}^{d} H_{a_{i},\alpha_{i}})}e^{-r+\sqrt{2r}z} \Big)\Big\}d\Phi(z).
\end{eqnarray*}

\end{theorem}

\begin{theorem}
\label{th3}
Let $X(\t)$ be a centered homogenous Gaussian field with \pzx{unit variance} and covariance function $r(\t)$ satisfying $\mathbf{A1}-\mathbf{A3}$. Then for any dense grids $ \Re(\delta_{i}) $, $ i=1,2,\cdots,d $, we have
\begin{eqnarray*}
&&\lim_{\T\to\infty}\P\left(\mathbf M_{\T}\leq\mathbf u_{\T},\mathbf m_{\T}\leq\mathbf v_{\T}\right)\\
&=&\int_{-\infty}^{+\infty}\Big[1-\exp\left(-e^{x_{1}+r+\sqrt{2r}z}\right)-\exp\left(-e^{y_{1}+r+\sqrt{2r}z}\right)
+\exp\left(-e^{\max(x_{1},y_{1})+r+\sqrt{2r}z}\right)\Big]\\
&&\times\exp\Big\{-\left(e^{-\min(x_{2},y_{2})-r+\sqrt{2r}z}\right)\Big\}d\Phi(z),
\end{eqnarray*}
where  $\mathbf{u}_{\T}$, $\mathbf{v}_{\T}$  are given by \eqref{eq4.2}.
\end{theorem}

\begin{remark}
Similar to the weakly dependent stationary Gaussian sequences and processes, for homogeneous Gaussian random fields, Theorems \ref{th1}-\ref{th3} shows that $\mathbf M_{\T}$ and $\mathbf m_{\T}$ are asymptotically independent if $r=0$.
\end{remark}

\section{Auxiliary results}
\label{sec3}

For simplicity, let $u=a_{\T}=\sqrt{2\log(\prod_{i=1}^{d}T_{i})} $, so for $ i=1,2,\cdots,d$, if $\delta_{i}u^{2/\alpha_{i}} \to 0$,  \pzx{the grid is dense};  if $\delta_{i}u^{2/\alpha_{i}} \to \infty$ , the grid is sparse, and $\delta_{i}u^{2/\alpha_{i}} \to D_{i}\in (0,\infty)$ for the Pickands grid.

For the sparse grids, let $ \delta_{i}=\delta_{i}(u)=l_{i}(u)u^{-2/\alpha_{i}}, i=1,2,\cdots,d $, where $ l_{i}(u)\to\infty $ as $ u\to\infty $ with $ \delta_{i}(u)\leq \delta_{0} $ for some $ \delta_{0}>0 $. Denote $ \I_{(\boldsymbol{\delta})}=\prod_{i=1}^{d}[-\delta_{i},\delta_{i}] $ and
\begin{equation}\label{eqpux}
P(u,x)=\P\left( X(\mathbf{0})>u, \max_{\t\in\I_{(\boldsymbol{\delta})}} X(\t)> u+\frac{v^2+x}{u}\right),
\end{equation}
where $ v^{2}= (\sum_{i=1}^{d} {2}/{\alpha_{i}})\log u + \log (\prod_{i=1}^{d} \delta_{i}).$ The following Lemmas \ref{lem1}-\ref{lem3} extended Lemmas 1, 2 and 4 of Piterbarg (2004) from stationary Gaussian processes to homogenous Gaussian random fields.

\begin{lemma}
\label{lem1}
Suppose $ \Re(\delta_{i})$, $ i=1,2,\cdots,d $ are all sparse grids and the conditions $ \mathbf{A1}-\mathbf{A2} $ hold, then $P(u,x)$ given by \eqref{eqpux} satisfies $P(u,x)= o(\Psi(u))$ as $ u\to\infty $.
\end{lemma}

\pzx{\noindent {\bf Proof.} The proof is similar to Lemma 1 of Piterbarg (2004) and Lemma A1 of Tan and Wang (2015).}
\qed

Now, define
\[ P_{\S}(u,x)=\P\left( \max_{\t\in\I_{\S}\cap\prod_{i=1}^{d} \Re(\delta_{i})} X(\t)> u, \max_{\t\in\I_{\S}} X(\t)> u+\frac{v^2+x}{u}\right)\]
and for small $ \epsilon>0 $
\begin{eqnarray}
\label{eq3.1}
\delta(\epsilon)=\inf_{\max\{|t_{i}|,i=1,2,\cdots,d\}>\epsilon} (1-r(\t)) >0.
\end{eqnarray}

\begin{lemma}
\label{lem2}
Suppose $ \Re(\delta_{i})$, $ i=1,2,\cdots,d,$ are all sparse grids and the conditions $ \mathbf{A1}-\mathbf{A2} $ hold. Let $ S_{i}=S_{i}(u)\geq 2\delta_{i} $, $ i=1,2,\cdots,d $, for all $u$, and $ \prod_{i=1}^{d}S_{i}u^{\frac{2}{\alpha_{i}}}=o\left(\exp(u^{2}\delta(\epsilon))\right) $ as $ u\to\infty $. Then there exists an $ \epsilon>0 $ such that
\begin{eqnarray}
\label{eq3.2}
\P\left(\max_{\t\in\I_{\S}} X(\t)> u+\frac{v^2+x}{u}\right) \sim \left(\prod_{i=1}^{d}S_{i}\delta_{i}^{-1}H_{\alpha_{i}}\right) e^{-x} \Psi(u),
\end{eqnarray}

\begin{eqnarray}
\label{eq3.3}
\P\left( \max_{\t\in\I_{\S}\cap\prod_{i=1}^{d} \Re(\delta_{i})} X(\t)> u\right)\sim \left(\prod_{i=1}^{d}S_{i}\delta_{i}^{-1}\right) \Psi(u)
\end{eqnarray}
as $ u\to\infty $ and
\[ P_{\S}(u,x)=o\left(\P\left( \max_{\t\in\I_{\S}\cap\prod_{i=1}^{d} \Re(\delta_{i})} X(\t)> u\right)+\P\left( \max_{\t\in\I_{\S}} X(\t)> u+\frac{v^2+x}{u}\right)\right)\]
as $ u\to\infty $, so that
\begin{eqnarray}
\label{eq3.4}
&& 1-\P\left( \max_{\t\in\I_{\S}\cap\prod_{i=1}^{d} \Re(\delta_{i})} X(\t)\leq u, \max_{\t\in\I_{\S}} X(\t)\leq u+\frac{v^2+x}{u}\right)\nonumber\\
&&\sim\P\left( \max_{\t\in\I_{\S}\cap\prod_{i=1}^{d} \Re(\delta_{i})} X(\t)> u\right)+\P\left( \max_{\t\in\I_{\S}} X(\t)> u+\frac{v^2+x}{u}\right)
\end{eqnarray}
as $ u\to\infty $.
\end{lemma}

\pzx{\noindent {\bf Proof.} The proof is similar to Lemma 2 of Piterbarg (2004) and Lemma A2 of Tan and Wang (2015).}
\qed

In the following Lemma, we can prove that the maxima asymptotically coincide when the grids all are  dense.
\begin{lemma}
\label{lem3}
Let $ S_{i}=S_{i}(u) $ with $ (\prod_{i=1}^{d} S_{i}u^{-2/\alpha_{i}})\to\infty $ and $ \prod_{i=1}^{d} S_{i}=O(\exp(\kappa u^{2})) $ with $ \kappa\in (0,1/2) $ as $ u\to\infty $. For any dense grids $ \Re(\delta_{i})=\{ a_{i}ku^{-2/\alpha_{i}},\;\  k\in N\} $, $ i=1, 2, \cdots, d$ , we have
\begin{eqnarray}
\label{eq4.1}
\P\left( \max_{\t\in\I_{\S}\cap\prod_{i=1}^{d} \Re(\delta_{i})} X(\t)\leq u\right)-\P\left( \max_{\t\in\I_{\S}} X(\t)\leq u\right)\leq \rho_{\a} \left(\prod_{i=1}^{d} H_{\alpha_{i}} S_{i} u^{\frac{2}{\alpha_{i}}}\right) \Psi(u),
\end{eqnarray}
where $ \rho_{\a}\to 0 $ as $ \a\to 0$.
\end{lemma}
\noindent {\bf Proof.} This lemma follows from the  Lemma D.1 and  Lemma 15.3 of Piterbarg (1996), we can also find the detailed proof in the Lemma 12.3.2 of Leadbetter (1983).
\qed

\label{sec3}
Following Tan and Wang (2015), define $ \rho(\T)=r/ \log(\prod_{i=1}^{d} T_{i}) $ and let $ 1>a>b>0 $ be constants. Dividing $ [0,T_{i}] $ into intervals with length $ T_{i}^{a} $ alternating with shorter intervals with length $ T_{i}^{b}$, $ i=1,2,\cdots,d $. Then the number of the long intervals is at most \pzx{$ n_{i}=\lfloor T_{i}/(T_{i}^{a}+T_{i}^{b})\rfloor $. Here $ \lfloor \cdot \rfloor $ represents the integer part of the real number}. Denote $ \O_{\i}=\prod_{j=1}^{d}[(i_{j}-1)(T_{j}^{a}+T_{j}^{b}), (i_{j}-1)(T_{j}^{a}+T_{j}^{b})+T_{j}^{a}] $, $ \E_{\i}=\prod_{j=1}^{d}[(i_{j}-1)(T_{j}^{a}+T_{j}^{b}), i_{j}(T_{j}^{a}+T_{j}^{b})] $, $ \i=\mathbf{1},\cdots,\mathbf{n} $ and $\O=\bigcup_{\i}\O_{\i}$. Let $ \{ X_{\i}(\t),\t\geq \mathbf{0}\} $, $ \i \geq 1 $ be independent copies of $ \{X(\t),\t\geq \mathbf{0}\} $ and $\{ \eta(\t),\t\geq \mathbf{0}\} $ be such that $ \eta(\t)=X_{\i}(\t) $ for $ \t\in \E_{\i} $.

Define
\[ \xi_{\T}(\t)=(1-\rho(\T))^{1/2}\eta(\t)+\rho^{1/2}(\T)U,\;\ \t \in \mathbf{I}_{\T},\]
where $ U $ is a standard normal variable independent of $ \{ \eta(\t),\t\geq 0\} $. Then $ \varrho(\s,\t) $, covariance function of $ \{\xi_{\T}(\t),\;\ \t \in \I_{\T}\} $,  is
\begin{eqnarray*}
\varrho(\s,\t) &=& \left\{ {{\begin{array}{*{20}c}
 {r(\t,\s)+\left(1-r(\t,\s)\right)\rho(\T),\;\ \s\in\E_{\i},\;\ \t\in\E_{\j},\;\ \i=\mathbf{j};} \hfill \\
 {\rho(\T),\quad\quad\quad\quad\qquad\quad\quad\quad\quad \s\in\E_{\i},\;\ \t\in\E_{\j},\;\ \i\neq\j. } \hfill \\
\end{array} }} \right.\\
\end{eqnarray*}
Let
\[\mathbf{ M}_{\O}=\left( M_{\O},M_{\O}^{\boldsymbol{\delta}}\right)=\left( \max_{\t\in \O} X(\t),\max_{\t\in \O \bigcap \prod_{i=1}^{d} \Re(\delta_{i})} X(\t) )\right) ,\]
\[\mathbf{ m}_{\O}=\left( m_{\O},m_{\O}^{\boldsymbol{\delta}}\right)=\left( \min_{\t\in \O} X(\t),\min_{\t\in \O \bigcap \prod_{i=1}^{d} \Re(\delta_{i})} X(\t) )\right) .\]

\begin{lemma}
\label{lem4}
Suppose $ \Re(\delta_{i}), i=1,2,\cdots,d $ are all sparse grids or all Pickands grids. Then for the $\mathbf{u}_{\T}$, $\mathbf{v}_{\T}$ given by \eqref{eq4.2}  we have
\begin{eqnarray*}
 |\P\left( \mathbf{v}_{\T}<\mathbf{m}_{\T}\leq\mathbf{M}_{\T}\leq\mathbf{u}_{\T}\right)-\P\left( \mathbf{v}_{\T}<\mathbf{m}_{\O}\leq\mathbf{M}_{\O}\leq\mathbf{u}_{\T}\right)| \to 0
 \end{eqnarray*}
as $ \T\to\infty $.
\end{lemma}

\noindent {\bf Proof.}  Note that the homogeneous Gaussian fields $ \{X(\t),\t\geq \mathbf{0}\} $ and $ \{-X(\t),\t\geq \mathbf{0}\} $ have the same distribution so that
\begin{eqnarray*}
&&|\P\left( \mathbf{v}_{\T}<\mathbf{m}_{\T}\leq\mathbf{M}_{\T}\leq\mathbf{u}_{\T}\right)-\P\left( \mathbf{v}_{\T}<\mathbf{m}_{\O}\leq\mathbf{M}_{\O}\leq\mathbf{u}_{\T}\right)| \\
&&\leq \P\left( \max_{\t\in\I_{\T}\backslash \mathbf{O}} X(\t)> u_{\T}(x_{2})\right)+ \P\left( \max_{\t\in(\I_{\T}\backslash \mathbf{O})\cap\prod_{i=1}^{d} \Re(\delta_{i})} X(\t)> u_{\T}^{\boldsymbol{\delta}}(y_{2})\right)\\
&&+\P\left( \max_{\t\in\I_{\T}\backslash \mathbf{O}} (-X(\t))> u_{\T}(-x_{1})\right)+ \P\left( \max_{\t\in(\I_{\T}\backslash \mathbf{O})\cap\prod_{i=1}^{d} \Re(\delta_{i})} (-X(\t))> u_{\T}^{\boldsymbol{\delta}}(-y_{1})\right).
\end{eqnarray*}
By arguments similar to  Lemma 3.1 of Tan and Wang (2015)  (denote by $ mes(\cdot) $ the Lebesgue measure), we have
\begin{eqnarray*}
 \P\left( \max_{\t\in\I_{\T}\backslash \mathbf{O}} X(\t)> u_{\T}(x_{2})\right)= O(1) mes(\mathbf{I}_{\T}\backslash \mathbf{O})(( u_{\T}(x_{2}))^{\sum_{i=1}^{d} \frac{2}{\alpha_{i}}}\Psi( u_{\T}(x_{2}))
 \to 0
 \end{eqnarray*}
as $\T\to\infty$. Combining \eqref{eq3.3}, we can get the assertion of this lemma.\qed

\begin{lemma}
\label{lem5}
Denote the grids $ \Re(q_{i})$ with $ q_{i}=\gamma_{i}u^{-\frac{2}{\alpha_{i}}} $, $ \gamma_{i}>0 $ and $ i=1,2,\cdots,d $. Assume that $\Re(\delta_{i}), i=1,2,\cdots,d,$ are all sparse grids or all Pickands grids. Then for  ${\u}_{\T}$ and $\v_{\T}$ given by \eqref{eq4.2} we have
\begin{eqnarray}
\label{eq5.1}
&&\Big| \P\Big( v_{\T}(x_{1}) < \min_{\t\in \O} X(\t) \leq \max_{\t\in \O} X(\t)\leq u_{\T}(x_{2}), \nonumber\\ && \qquad v_{\T}^{\boldsymbol{\delta}}(y_{1}) < \min_{\t\in \O \cap\prod_{i=1}^{d} \Re(\delta_{i})} X(\t) \leq \max_{\t\in \O\cap\prod_{i=1}^{d} \Re(\delta_{i})} X(\t)\leq  u_{\T}^{\boldsymbol{\delta}}(y_{2}) \Big) \nonumber \\
&&\qquad\qquad-\P\Big( v_{\T}(x_{1}) < \min_{\t\in \O\cap\prod_{i=1}^{d} \Re(q_{i})} X(\t) \leq \max_{\t\in \O\cap \prod_{i=1}^{d} \Re(q_{i})} X(\t)\leq u_{\T}(x_{2}), \nonumber\\ &&\qquad\qquad\qquad v_{\T}^{\boldsymbol{\delta}}(y_{1}) < \min_{\t\in \O \cap\prod_{i=1}^{d} \Re(\delta_{i})} X(\t) \leq \max_{\t\in \O\cap\prod_{i=1}^{d} \Re(\delta_{i})} X(\t)\leq  u_{\T}^{\boldsymbol{\delta}}(y_{2})\Big) \Big| \nonumber \\
&&\to 0
\end{eqnarray}
as $ \T\to\infty $ and $ \boldsymbol{\gamma}=(\gamma_{1},\gamma_{2},\cdots,\gamma_{d})\to 0 $.
\end{lemma}

\noindent {\bf Proof.} It follows from Lemma \ref{lem3} and the fact $ \{X(\t),\t\geq \mathbf{0}\} =^{d} \{-X(\t),\t\geq \mathbf{0}\} $ that \pzx{the left hand side of \eqref{eq5.1}} can be bounded by
\begin{eqnarray*}
&& \P\left( \max_{\t\in \O \cap\prod_{i=1}^{d} \Re(q_{i})} X(\t)\leq u_{\T}(x_{2}), \max_{\t\in \O} X(\t) > u_{\T}(x_{2}) \right)\\
&&+ \P\left( \max_{\t\in \O\cap \prod_{i=1}^{d} \Re(q_{i})} (-X(\t))> u_{\T}(-x_{1}), \max_{\t\in \O} (-X(\t)) < u_{\T}(-x_{1}) \right)\\
&\to& 0
\end{eqnarray*}
as $ \T\to\infty $ and $ \boldsymbol{\gamma}=(\gamma_{1},\gamma_{2},\cdots,\gamma_{d})\to 0 $. The result follows.
\qed

\begin{lemma}
\label{lem6}
Suppose $ \Re(\delta_{i}), i=1,2,\cdots,d $ are all sparse grids or all Pickands grids, then for  ${\u}_{\T}$ and ${\v}_{\T}$ given by \eqref{eq4.2} we have
\begin{eqnarray}
\label{eq6.1}
&&\Big|\P\Big( v_{\T}(x_{1}) < \min_{\t\in \O\cap\prod_{i=1}^{d} \Re(q_{i})} X(\t) \leq \max_{\t\in \O\cap \prod_{i=1}^{d} \Re(q_{i})} X(\t)\leq u_{\T}(x_{2}), \nonumber\\
&&\qquad v_{\T}^{\boldsymbol{\delta}}(y_{1}) < \min_{\t\in \O \cap\prod_{i=1}^{d} \Re(\delta_{i})} X(\t) \leq \max_{\t\in \O\cap\prod_{i=1}^{d} \Re(\delta_{i})} X(\t)\leq  u_{\T}^{\boldsymbol{\delta}}(y_{2})\Big) \nonumber\\
&&\qquad\qquad-\P\Big( v_{\T}(x_{1}) < \min_{\t\in \O\cap\prod_{i=1}^{d} \Re(q_{i})} \xi_{\T}(\t) \leq \max_{\t\in \O \cap\prod_{i=1}^{d} \Re(q_{i})} \xi_{\T}(\t)\leq u_{\T}(x_{2}), \nonumber\\
&&\qquad\qquad\qquad v_{\T}^{\boldsymbol{\delta}}(y_{1}) < \min_{\t\in \O \cap\prod_{i=1}^{d} \Re(\delta_{i})} \xi_{\T}(\t) \leq \max_{\t\in \O\cap\prod_{i=1}^{d} \Re(\delta_{i})} \xi_{\T}(\t)\leq  u_{\T}^{\boldsymbol{\delta}}(y_{2})) \Big| \nonumber\\
&&\to 0
\end{eqnarray}
as $ \T\to\infty $.
\end{lemma}

\noindent {\bf Proof.} It follows from the Normal Comparison Lemma \pzx{(see e.g. Leadbetter et al. (1983) and Li and Shao (2002)) } that the left hand side of \eqref{eq6.1} can be bounded by
\begin{eqnarray*}
&& C\sum_{\mbox{\tiny$\begin{array}{c} \k\q\in \O_{\i}, \l\q\in \O_{\j} \\ \k\q\neq \l\q, \mathbf{1}\leq \i,\j\leq \mathbf{n}\end{array}$}} | r(\k\q,\l\q)-\varrho(\k\q,\l\q)|\int_{0}^{1}\frac{1}{\sqrt{1-r^{h}(\k\q,\l\q)}}\times \\
&& \Big[\exp\left(-\frac{u_{\T}^{2}(x_{2})}{1+r^{h}(\k\q,\l\q)}\right)
+\exp\left(-\frac{v_{\T}^{2}(x_{1})}{1+r^{h}(\k\q,\l\q)}\right)\Big] dh \\
&& +C\sum_{\mbox{\tiny$\begin{array}{c} \k\boldsymbol{\delta}\in \O_{\i}, \l\boldsymbol{\delta}\in \O_{\j} \\ \k\boldsymbol{\delta}\neq \l\boldsymbol{\delta}, \mathbf{1}\leq \i,\j\leq \mathbf{n}\end{array}$}} | r(\k\boldsymbol{\delta},\l\boldsymbol{\delta})-\varrho(\k\boldsymbol{\delta},\l\boldsymbol{\delta})|\int_{0}^{1}\frac{1}{\sqrt{1-r^{h}(\k\boldsymbol{\delta},\l\boldsymbol{\delta})}}\times \\
&& \Big[\exp\left(-\frac{\left(u_{\T}^{\boldsymbol{\delta}}(y_{2})\right)^{2}}{1+r^{h}(\k\boldsymbol{\delta},\l\boldsymbol{\delta})}\right)
+\exp\left(-\frac{\left(v_{\T}^{\boldsymbol{\delta}}(y_{1})\right)^{2}}{1+r^{h}(\k\boldsymbol{\delta},\l\boldsymbol{\delta})}\right)\Big] dh \\
 && +C\sum_{\mbox{\tiny$\begin{array}{c} \k\q\in \O_{\i}, \l\boldsymbol{\delta}\in \O_{\j} \\ \k\q\neq \l\boldsymbol{\delta}, \mathbf{1}\leq \i,\j\leq \mathbf{n}\end{array}$}} | r(\k\q,\l\boldsymbol{\delta})-\varrho(\k\q,\l\boldsymbol{\delta})|\int_{0}^{1}\frac{1}{\sqrt{1-r^{h}(\k\q,\l\boldsymbol{\delta})}}\times \\
&& \Big[\exp\left(-\frac{u_{\T}^{2}(x_{2})+\left(u_{\T}^{\boldsymbol{\delta}}(y_{2})\right)^{2}}{2\left(1+r^{h}(\k\q,\l\boldsymbol{\delta})\right)}\right)
+\exp\left(-\frac{v_{\T}^{2}(x_{1})+\left(v_{\T}^{\boldsymbol{\delta}}(y_{1})\right)^{2}}{2\left(1+r^{h}(\k\q,\l\boldsymbol{\delta})\right)}\right)\\
&&+\exp\left(-\frac{v_{\T}^{2}(x_{1})+\left(u_{\T}^{\boldsymbol{\delta}}(y_{2})\right)^{2}}{2\left(1+r^{h}(\k\q,\l\boldsymbol{\delta})\right)}\right)
+\exp\left(-\frac{u_{\T}^{2}(x_{2})+\left(v_{\T}^{\boldsymbol{\delta}}(y_{1})\right)^{2}}{2\left(1+r^{h}(\k\q,\l\boldsymbol{\delta})\right)}\right)\Big] dh\\
&&:= F_{1}+F_{2}+F_{3},
\end{eqnarray*}
where \pzx{$ r^{h}(\mathbf{x},\mathbf{y})=h r(\mathbf{x},\mathbf{y})+(1-h)\varrho(\mathbf{x},\mathbf{y})  $ }with $ h\in [0,1]$.

Next we will show that $ F_{1}\to 0 $, $ F_{2}\to 0 $ and $ F_{3}\to\ 0 $ as $ \T\to \infty $, respectively. For $ F_{1} $, we first consider the case that $ \k\q $, $ \l\q $ in the same interval firstly, and split $ F_{1} $ into the following two parts:
\begin{eqnarray}
\label{eq6.2}
\sum_{\mbox{\tiny$\begin{array}{c} \k\q, \l\q \in \O_{\i}, \k\q\neq\l\q, \i=\mathbf{1,2,\cdots,n} \\ \max\{|l_{j}q_{j}-k_{j}q_{j}|,j=1,2,\cdots,d\}\leq \epsilon\end{array}$}} + \sum_{\mbox{\tiny$\begin{array}{c} \k\q, \l\q \in \O_{\i}, \k\q\neq\l\q, \i=\mathbf{1,2,\cdots,n} \\ \max\{|l_{j}q_{j}-k_{j}q_{j}|,j=1,2,\cdots,d\} > \epsilon\end{array}$}} := J_{\T,1} +J_{\T,2}.
\end{eqnarray}
Following condition $ \mathbf{A1} $, we can choose small enough $ \epsilon >0 $ such that $ \max\{|l_{j}q_{j}-k_{j}q_{j}|,j=1,2,\cdots,d\}\leq \epsilon $ and for all $  |t_{i}|\leq \epsilon < 2^{-1/\alpha_{i}} $,
\begin{eqnarray}
\label{eq6.3}
\frac{1}{2} \left(\sum_{i=1}^{d} |t_{i}|^{\alpha_{i}}\right)\leq 1-r(\t) \leq 2 \left(\sum_{i=1}^{d} |t_{i}|^{\alpha_{i}}\right),
\end{eqnarray}
then, by definition of $ \xi_{\T}(\t) $, we have $ \varrho(\k\q,\l\q)-r(\k\q,\l\q)=\rho(\T)\left( 1-r(\k\q,\l\q)\right) $ and $ \varrho(\k\q,\l\q) \sim r(\k\q,\l\q) $ for sufficiently large $ \T $. With  $v_{\T}$ and $u_{\T}$ given by \eqref{eq4.2} we have
\begin{align}
\label{eq6.4}
v_{\T}^{2}(x_{1})=u^{2}-\left[1-\left(\sum_{i=1}^{d} \frac{2}{\alpha_{i}}\right)\right]\log\left(\frac{u^{2}}{2}\right)+O(1),\quad
u_{\T}^{2}(x_{2})=u^{2}-\left[1-\left(\sum_{i=1}^{d} \frac{2}{\alpha_{i}}\right)\right]\log\left(\frac{u^{2}}{2}\right)+O(1).
\end{align}
Hence,
\begin{eqnarray}
\label{eq6.5}
J_{\T,1}&\leq& C \sum_{\mbox{\tiny$\begin{array}{c} \k\q, \l\q \in \O_{\i}, \k\q\neq\l\q, \i=\mathbf{1,2,\cdots,n} \\ \max\{|l_{j}q_{j}-k_{j}q_{j}|,j=1,2,\cdots,d\}\leq \epsilon\end{array}$}} | r(\k\q,\l\q)-\varrho(\k\q,\l\q)| \frac{1}{\sqrt{1-r(\k\q,\l\q)}}
\exp\left(-\frac{u_{\T}^{2}(x_{2})}{1+r(\k\q,\l\q)}\right) \nonumber\\
&\leq & C (\prod_{i=1}^{d} T_{i}q_{i}^{-1}) \rho(\T) \sum_{j=1}^{d}  \sum_{ 0\leq k_{j}q_{j} \leq \epsilon} \sqrt{1-r(\k\q)}  \exp\left(-\frac{u_{\T}^{2}(x_{2})}{1+r(\k\q)}\right) \nonumber\\
&\leq & C (\prod_{i=1}^{d} T_{i}q_{i}^{-1}) \rho(\T)\exp(-\frac{u_{\T}^{2}(x_{2})}{2})\sum_{j=1}^{d}  \sum_{ 0\leq k_{j}q_{j} \leq \epsilon} \sqrt{1-r(\k\q)}  \exp\left(-\frac{(1-r(\k\q))u_{\T}^{2}(x_{2})}{2(1+r(\k\q))}\right) \nonumber\\
&\leq & C u^{-1} \sum_{j=1}^{d}  \sum_{ 0\leq k_{j}q_{j} \leq \epsilon} \left[\sum_{i=1}^{d} (k_{i}q_{i})^{\alpha_{i}}\right]^{1/2}\exp\left(-\frac{1}{4}\left[\sum_{i=1}^{d} (k_{i}q_{i})^{\alpha_{i}}\right]  \left(\log \prod_{i=1}^{d} T_{i}\right) \right) \nonumber\\
&\leq & C u^{-1} \sum_{j=1}^{d}  \sum_{ 0\leq k_{j}q_{j} \leq \epsilon} \exp\left(-\frac{1}{4}\left[\sum_{i=1}^{d} (k_{i}q_{i})^{\alpha_{i}}\right]  \left(\log \prod_{i=1}^{d} T_{i}\right) \right) \nonumber\\
&\leq & C u^{-1} \sum_{j=1}^{d} \sum_{k_{j}=0}^{\infty} e^{-\frac{1}{4}(k_{j}\gamma_{j})^{\alpha_{j}}} \nonumber\\
&\leq & C u^{-1}\to0
\end{eqnarray}
as $ u\to\infty $.

Let
$\varpi(\t,\s)=\max\{ |r(\t,\s)|,|\varrho(\t,\s)| \} $ and
\[\theta(z)=\sup_{\mbox{\tiny$\begin{array}{c} \mathbf{0}\leq \s,\t\leq \T \\ \max\{|s_{i}-t_{i}|,j=1,2,\cdots,d\}> z\end{array}$}}\{ \varpi(\t,\s) \}.\]
For $ J_{\T,2} $, by the fact that $ u_{\T}^{2}(x_{2})\sim 2 \log \prod_{i=1}^{d} T_{i}=u^{2}$ we have
\begin{eqnarray}
\label{eq6.6}
J_{\T,2}&\leq& C \sum_{\mbox{\tiny$\begin{array}{c} \k\q, \l\q \in \O_{\i}, \k\q\neq\l\q, \i=\mathbf{1,2,\cdots,n} \\ \max\{|l_{j}q_{j}-k_{j}q_{j}|,j=1,2,\cdots,d\}> \epsilon\end{array}$}} | r(\k\q,\l\q)-\varrho(\k\q,\l\q)| \exp\left(-\frac{u_{\T}^{2}(x_{2})}{1+\varpi(\k\q,\l\q)}\right)\nonumber\\
&\leq & C (\prod_{i=1}^{d} T_{i}q_{i}^{-1}) \sum_{j=1}^{d} \sum_{\mbox{\tiny$\begin{array}{c} 0\leq k_{j}q_{j}\leq T_{j}^{a} \\ \max\{k_{j}q_{j},j=1,2,\cdots,d\}> \epsilon\end{array}$}} \exp\left(-\frac{u_{\T}^{2}(x_{2})}{1+\theta(\epsilon)}\right)   \nonumber\\
&\leq & C \exp\left[\left(a-\frac{1-\theta(\epsilon)}{1+\theta(\epsilon)}\right)\frac{u^{2}}{2}+\sum_{i=1}^{d} \frac{4}{\alpha_{i}}\log u\right]\to0
\end{eqnarray}
as $ u\to\infty $ since $ a<\frac{1-\theta(\epsilon)}{1+\theta(\epsilon)} $. Combining with \eqref{eq6.2}, \eqref{eq6.5} and \eqref{eq6.6}, we show that $F_1\to0$ for the first case that $ \k\q $, $ \l\q $ in the same interval.

Second, we consider the case that $ \k\q\in \O_{\i} $, $ \l\q\in \O_{\j} $, $ \i\neq \j $. Note that the distance between the points in any two rectangles $ \O_{\i} $ and $ \O_{\j} $ are larger than $ \min\{ T_{i}^{b}, i=1,2,\cdots,d \} $ and $ \varrho(\k\q,\l\q)=\rho(\T) $ for $ \k\q\in \O_{\i} $, $ \l\q\in \O_{\j} $, $ \i\neq \j $. Then, $ F_{1} $ can be bounded by
\begin{eqnarray}
\label{eq6.7}
C \sum_{\mbox{\tiny$\begin{array}{c} \k\q\in \O_{\i},\l\q\in \O_{\j} \\ \mathbf{1\leq i\neq j\leq n}\end{array}$}} | r(\k\q,\l\q)-\rho(\T)| \exp\left(-\frac{u_{\T}^{2}(x_{2})}{1+\varpi(\k\q,\l\q)}\right).
\end{eqnarray}
Split \eqref{eq6.7} into two parts, the first for $ \min\{|k_{j}q_{j}-l_{j}q_{j}|,j=1,2,\cdots,d\}> 0 $, the second for $ \min\{|k_{j}q_{j}-l_{j}q_{j}|,j=1,2,\cdots,d\}= 0 $ and denote them $ S_{\T,1} $, $ S_{\T,2} $, respectively. Let $ \beta $ be such that $ 0<b<a<\beta< \frac{1-\theta(\epsilon)}{1+\theta(\epsilon)} $ for all sufficiently large $ \T $.

For $ S_{\T,1} $, we can also divide it into the following two parts:
\[ S_{\T,1}=C\sum_{\mbox{\tiny$\begin{array}{c} \k\q\in \O_{\i},\l\q\in \O_{\j},\mathbf{1\leq i\neq j\leq n} \\ \prod_{m=1}^{d}|k_{m}q_{m}-l_{m}q_{m}|\leq {\prod_{i=1}^{d} T_{i}}^{\beta}\end{array}$}}+C\sum_{\mbox{\tiny$\begin{array}{c} \k\q\in \O_{\i},\l\q\in \O_{\j},\mathbf{1\leq i\neq j\leq n} \\ \prod_{m=1}^{d}|k_{m}q_{m}-l_{m}q_{m}|> \prod_{i=1}^{d} T_{i}^{\beta}\end{array}$}} :=S_{\T,11}+S_{\T,12} ,\]
then, by the same arguments as used in \eqref{eq6.6}, we have
\begin{eqnarray}
\label{eq6.8}
S_{\T,11}&\leq & C (\prod_{i=1}^{d} T_{i}q_{i}^{-1}) \sum_{j=1}^{d} \sum_{\mbox{\tiny$\begin{array}{c} 0\leq k_{j}q_{j}\leq T_{j} \\ \prod_{m=1}^{d}k_{m}q_{m}\leq {\prod_{i=1}^{d} T_{i}}^{\beta} \end{array}$}} \exp\left(-\frac{u_{\T}^{2}(x_{2})}{1+\theta(\epsilon)}\right)   \nonumber\\
&\leq & C \exp\left[\left(\beta-\frac{1-\theta(\epsilon)}{1+\theta(\epsilon)}\right)\frac{u^{2}}{2}+\sum_{i=1}^{d} \frac{4}{\alpha_{i}}\log u\right],\nonumber\\
&\to& 0
\end{eqnarray}
as $u\to\infty  $ since $\beta< \frac{1-\theta(\epsilon)}{1+\theta(\epsilon)}$.

To deal with $ S_{\T,12} $, we can define
$ w_{1}(\t)=\max\{ |r(\t)|,|\rho(\T)| \}$ and
\[\theta_{1}(\mathbf{z})=\sup_{\mbox{\tiny$\begin{array}{c} \mathbf{0}\leq \t\leq \T \\ \prod_{i=1}^{d}|t_{i}|> \prod_{i=1}^{d} z_{i} \end{array}$}}\{w_{1}(\t)\}. \]
 Denote  $ \T^{\beta}=(T_{1}^{\beta},T_{1}^{\beta},\cdots,T_{d}^{\beta}) $ and by the condition $ \mathbf{A3} $, we have $ \theta_{1}(\T^{\beta}) \leq Cu^{-2}$ for sufficiently large $ \T $ and $ \prod_{j=1}^{d} k_{j}q_{j} >{\prod_{i=1}^{d} T_{i}}^{\beta}=e^{\beta u^{2}/2}$, then, using \eqref{eq6.4}, we have
\begin{eqnarray}
\label{eq6.9}
&&\frac{(\prod_{i=1}^{d} T_{i}q_{i}^{-1})^{2}}{\log\left(\prod_{i=1}^{d} T_{i}\right)}\exp\left(-\frac{u_{\T}^{2}(x_{2})}{1+\theta_{1}(\T^{\beta})}\right)\nonumber\\
&\le& O(1)\exp\left[\frac{C}{C/u^{2}+\beta/2}+2C\left(\sum_{i=1}^{d}2/\alpha_{i}-1\right)\frac{\log u}{C+\beta u^{2}/2}\right]=O(1)
\end{eqnarray}
as $ u\to\infty $. Therefore, we have
\begin{eqnarray}
\label{eq6.10}
S_{\T,12}&\leq& C (\prod_{i=1}^{d} T_{i}q_{i}^{-1}) \sum_{\mbox{\tiny$\begin{array}{c} \mathbf{0}\leq \k\q \leq \T \\ \prod_{m=1}^{d}k_{m}q_{m}> {\prod_{i=1}^{d} T_{i}}^{\beta} \end{array}$}} |r(\k\q)-\rho(\T) |\exp\left(-\frac{u_{\T}^{2}(x_{2})}{1+\theta_{1}(\T^{\beta})}\right)   \nonumber\\
&\leq & C\frac{(\prod_{i=1}^{d} T_{i}q_{i}^{-1})^{2}}{\log\left(\prod_{i=1}^{d} T_{i}\right)}\exp\left(-\frac{u_{\T}^{2}(x_{2})}{1+\theta_{1}(\T^{\beta})}\right) \frac{\log\left(\prod_{i=1}^{d} T_{i}\right)}{(\prod_{i=1}^{d} T_{i}q_{i}^{-1})} \sum_{\mbox{\tiny$\begin{array}{c} \mathbf{0}\leq \k\q \leq \T \\ \prod_{m=1}^{d}k_{m}q_{m}> {\prod_{i=1}^{d} T_{i}}^{\beta} \end{array}$}} |r(\k\q)-\rho(\T) | \nonumber\\
&\leq & C \left(\prod_{i=1}^{d} q_{i} T_{i}^{-1}\right) \sum_{\mbox{\tiny$\begin{array}{c} \mathbf{0}\leq \k\q \leq \T \\ \prod_{m=1}^{d}k_{m}q_{m}> {\prod_{i=1}^{d} T_{i}}^{\beta} \end{array}$}} \left|r(\k\q)\log\left(\prod_{j=1}^{d} k_{j}q_{j}\right)-r \right| \nonumber\\
&+&  C \left(\prod_{i=1}^{d} q_{i} T_{i}^{-1}\right) \sum_{\mbox{\tiny$\begin{array}{c} \mathbf{0}\leq \k\q \leq \T \\ \prod_{m=1}^{d}k_{m}q_{m}> {\prod_{i=1}^{d} T_{i}}^{\beta} \end{array}$}} \left| 1-\frac{\log\prod_{i=1}^{d} T_{i}}{\log\prod_{j=1}^{d} k_{j}q_{j}} \right|.
\end{eqnarray}
Note that
\begin{eqnarray}
\label{eq6.11}
&&\left(\prod_{i=1}^{d} q_{i} T_{i}^{-1}\right) \sum_{\mbox{\tiny$\begin{array}{c} \mathbf{0}\leq \k\q \leq \T \\ \prod_{m=1}^{d}k_{m}q_{m}> {\prod_{i=1}^{d} T_{i}}^{\beta} \end{array}$}} \left| 1-\frac{\log\prod_{i=}^{d} T_{i}}{\log\prod_{j=1}^{d} k_{j}q_{j}} \right| \nonumber\\
&\leq & \frac{\prod_{i=1}^{d} q_{i} T_{i}^{-1}}{\log\prod_{i=}^{d} T_{i}^{\beta}} \sum \left|\log\left( \frac{\prod_{i=1}^{d} T_{i}}{\prod_{j=1}^{d} k_{j}q_{j}}\right)\right| \nonumber\\
&=& O\left( {u^{-2}}\int_{0}^{1}\cdots\int_{0}^{1} |\log\prod_{i=1}^{d} x_{i}| dx_{1}\cdots dx_{d}\right)\nonumber\\
&\to& 0
\end{eqnarray}
as $ u\to\infty $. Combining condition $ \mathbf{A3} $, and \eqref{eq6.9}-\eqref{eq6.11}, we can get $ S_{\T,1}\to 0 $ as $ \T\to\infty $.

Now, for $ S_{\T,2} $, we only prove the case that $ k_{1}q_{1}= l_{1}q_{1} $, and $ \min\{|k_{j}q_{j}-l_{j}q_{j}|,j=2,\cdots,d\}> 0  $. Other cases can be proved by the similar arguments. If $ \prod_{i=2}^{d} T_{i}\leq \prod_{i=1}^{d} T_{i}^{\beta}=e^{\beta u^{2}/2} $, where $ \beta< \frac{1-\theta(\epsilon)}{1+\theta(\epsilon)} $, by \pzx{the same arguments as used in \eqref{eq6.8} that}
\begin{eqnarray*}
S_{\T,2}&\leq & C (\prod_{i=1}^{d} T_{i}q_{i}^{-1})\sum_{j=2}^{d} \sum_{\mbox{\tiny$\begin{array}{c} 0\leq k_{j}q_{j}\leq T_{j} \\ k_{1}q_{1}=0\end{array}$}} |r(\k\q)-\rho(\T) |\exp\left(-\frac{u_{\T}^{2}(x_{2})}{1+\theta(\epsilon)}\right)   \nonumber\\
&\leq & C \exp\left[\left(\beta- \frac{1-\theta(\epsilon)}{1+\theta(\epsilon)}\right)u^{2}/2+ \left(\frac{1}{\alpha_{1}}+\sum_{i=2}^{d} \frac{2}{\alpha_{i}}\right)\log (u^{2}/2)\right] \nonumber\\
&\to& 0
\end{eqnarray*}
as $u\to\infty$. If $ \prod_{i=2}^{d} T_{i}> \prod_{i=1}^{d} T_{i}^{\beta}=e^{\beta u^{2}/2} $, split $ S_{\T,2} $ as follows:
\[ S_{\T,2}=C\sum_{\mbox{\tiny$\begin{array}{c} \k\q\in \O_{\i},\l\q\in \O_{\j},\mathbf{1\leq i\neq j\leq n} \\ \prod_{m=2}^{d}|k_{m}q_{m}-l_{m}q_{m}|\leq {\prod_{i=1}^{d} T_{i}}^{\beta}\\ k_{1}=l_{1}\end{array}$}}+C\sum_{\mbox{\tiny$\begin{array}{c} \k\q\in \O_{\i},\l\q\in \O_{\j},\mathbf{1\leq i\neq j\leq n} \\ \prod_{m=2}^{d}|k_{m}q_{m}-l_{m}q_{m}|> \prod_{i=1}^{d} T_{i}^{\beta}\\ k_{1}=l_{1}\end{array}$}} :=S_{\T,21}+S_{\T,22}.\]
For $ S_{\T,21} $, similar to \pzx{the arguments as used in \eqref{eq6.8}}, we have
\begin{eqnarray*}
S_{\T,21}&\leq & C (\prod_{i=1}^{d} T_{i}q_{i}^{-1}) \sum_{\mbox{\tiny$\begin{array}{c} \prod_{j=2}^{d} k_{j}q_{j}\leq \prod_{i=1}^{d}T_{i}^{\beta} \\ k_{1}q_{1}=0\end{array}$}} |r(\k\q)-\rho(\T) |\exp\left(-\frac{u_{\T}^{2}(x_{2})}{1+\theta(\epsilon)}\right)   \nonumber\\
&\leq & C \exp\left[\left(\beta- \frac{1-\theta(\epsilon)}{1+\theta(\epsilon)}\right)u^{2}/2+\left(\frac{1}{\alpha_{1}}+\sum_{i=2}^{d} \frac{2}{\alpha_{i}}\right)\log u^{2}/2\right] \nonumber\\
&\to& 0
\end{eqnarray*}
as $ u\to\infty $. To deal with $ S_{\T,22} $, let
$ w_{2}(\t)=\max\{ |r(0,t_{2},\cdots,t_{d})|,|\rho(\T)| \}$
and
\[\theta_{2}(\mathbf{z})=\sup_{\mbox{\tiny$\begin{array}{c} \mathbf{0}\leq \t\leq \T \\ \prod_{i=2}^{d}|t_{i}|> \prod_{i=1}^{d} z_{i} \end{array}$}}\{w_{2}(\t)\}. \]
By condition $ \mathbf{A3} $, we have $ \theta_{2}(\T^{\beta}) \leq C u^{-2} $ for sufficiently large $\T $ and $ \prod_{j=2}^{d} k_{j}q_{j} > \prod_{i=1}^{d} T_{i}^{\beta}=e^{\beta u^{2}/2} $, then, by \pzx{the same arguments as used in \eqref{eq6.9}},
\begin{eqnarray*}
\frac{(\prod_{i=1}^{d} T_{i}q_{i}^{-1})^{2}}{\log\left(\prod_{i=1}^{d} T_{i}\right)}\exp\left(-\frac{u_{\T}^{2}(x_{2})}{1+\theta_{2}(\T^{\beta})}\right)=O(1),
\end{eqnarray*}
thus we have
\begin{eqnarray}
\label{eq1}
S_{\T,22}&\leq & C (\prod_{i=1}^{d} q_{i}T_{i}^{-1})\log\left(\prod_{i=1}^{d} T_{i}\right) \sum_{\mbox{\tiny$\begin{array}{c} \prod_{j=2}^{d} k_{j}q_{j}> \prod_{i=1}^{d}T_{i}^{\beta} \\ k_{1}q_{1}=0\end{array}$}} \left| r(0,k_{2}q_{2},\cdots,k_{d}q_{d})-\rho(\T)  \right|\nonumber\\
&\leq & C (\prod_{i=1}^{d} q_{i}T_{i}^{-1})\log\left(\prod_{i=1}^{d} T_{i}\right) \sum_{\prod_{i=2}^{d} T_{i}\geq\prod_{j=2}^{d} k_{j}q_{j}> \prod_{i=1}^{d}T_{i}^{\beta}} \left| r(0,k_{2}q_{2},\cdots,k_{d}q_{d})+\rho(\T)  \right| \nonumber\\
&\leq& C\frac{q_{1}}{T_{1}},
\end{eqnarray}
which implies $ S_{\T,22}\to 0 $ as $ \T\to\infty $. Therefore, we showed that $ S_{\T,2}\to 0 $ as $ \T\to\infty $. Combining with $ S_{\T,1}\to0$, we showed that $ F_{1}\to 0 $ as $ \T\to\infty $ for the second case.

Arguments similar to the proof of $F_{1}\to 0$, we can show that $ F_{2}\to 0 $ as $ \T\to\infty $. Details are omitted here. The reminder is to show that $ F_{3}\to 0 $.

If $ \k\q $, $ \l\boldsymbol{\delta} $ in the same interval $ \O_{\i} $, split $ F_{3} $ into two parts as
\begin{eqnarray}
\label{eq6.12}
\sum_{\mbox{\tiny$\begin{array}{c} \k\q, \l\boldsymbol{\delta}\in \O_{\i}, \k\q\neq\l\boldsymbol{\delta}, \i=\mathbf{1,2,\cdots,n} \\ \max\{|l_{j}\delta_{j}-k_{j}q_{j}|,j=1,2,\cdots,d\}\leq \epsilon\end{array}$}} + \sum_{\mbox{\tiny$\begin{array}{c} \k\q, \l\boldsymbol{\delta} \in \O_{\i}, \k\q\neq\l\boldsymbol{\delta}, \i=\mathbf{1,2,\cdots,n} \\ \max\{|l_{j}\delta_{j}-k_{j}q_{j}|,j=1,2,\cdots,d\} > \epsilon\end{array}$}} := W_{\T,1} +W_{\T,2}.
\end{eqnarray}
Following condition $\mathbf{A1}$, we can choose small enough $ \epsilon >0 $ \pzx{so \eqref{eq6.3} is satisfied} and $ \max\{|l_{j}\delta_{j}-k_{j}q_{j}|,j=1,2,\cdots,d\}\leq \epsilon $. By definition of $ \xi_{\T}(\t) $, we have $ \varrho(\k\q,\l\boldsymbol{\delta})=r(\k\q,\l\boldsymbol{\delta})+\rho(\T)\left( 1-r(\k\q,\l\boldsymbol{\delta})\right) \sim r(\k\q,\l\boldsymbol{\delta}) $ for sufficient large $ \T $. If $ \Re(\delta_{i}), i=1,2,\cdots,d $ are all Pickands grids, by the arguments similar to the proof of $ F_{1} $, we can show $F_{3}\to0$.  If $ \Re(\delta_{i}), i=1,2,\cdots,d $  are all sparse grids, it follows from \eqref{eq4.2} that
\begin{eqnarray*}
 u_{\T}^{2}:&=&\frac{1}{2}\left(u_{\T}^{2}(x_{2})+(u_{\T}^{\boldsymbol{\delta}}(y_{2}))^{2}\right)=u^{2}-\left(1-\sum_{i=1}^{d} \frac{1}{\alpha_{i}}\right)\log (u^{2}/2)+ \log(\prod_{i=1}^{d} \delta_{i}^{-1})+O(1).
 \end{eqnarray*}
Similarly for $ \frac{1}{2}\left(v_{\T}^{2}(x_{1})+(v_{\T}^{\boldsymbol{\delta}}(y_{1}))^{2}\right) $, $ \frac{1}{2}\left(u_{\T}^{2}(x_{2})+(v_{\T}^{\boldsymbol{\delta}}(y_{1}))^{2}\right) $ and $ \frac{1}{2}\left(v_{\T}^{2}(x_{1})+(u_{\T}^{\boldsymbol{\delta}}(y_{2}))^{2}\right) $. In view of \eqref{eq6.3}, we have
 \begin{eqnarray}
 \label{eq6.13}
 W_{\T,1}&\leq&C \sum_{\mbox{\tiny$\begin{array}{c} \k\q, \l\boldsymbol{\delta} \in \O_{\i}, \k\q\neq\l\boldsymbol{\delta}, \i=\mathbf{1,2,\cdots,n} \\ \max\{|l_{j}\delta_{j}-k_{j}q_{j}|,j=1,2,\cdots,d\}\leq \epsilon\end{array}$}}  |r(\k\q,\l\boldsymbol{\delta})-\varrho(\k\q,\l\boldsymbol{\delta})| \frac{1}{\sqrt{1-r(\k\q,\l\boldsymbol{\delta})}}
\exp\left(-\frac{u_{\T}^{2}}{1+r(\k\q,\l\boldsymbol{\delta})}\right) \nonumber\\
&\leq & C  \rho(\T)   \exp\left(-\frac{u_{\T}^{2}}{2}\right)
 \sum_{\mbox{\tiny$\begin{array}{c} \k\q, \l\boldsymbol{\delta} \in \O_{\i}, \k\q\neq\l\boldsymbol{\delta}, \i=\mathbf{1,2,\cdots,n} \\ \max\{|l_{j}\delta_{j}-k_{j}q_{j}|,j=1,2,\cdots,d\}\leq \epsilon\end{array}$}} \sqrt{1-r(\k\q,\l\boldsymbol{\delta})}  \exp\left(-\frac{(1-r(\k\q,\l\boldsymbol{\delta}))u_{\T}^{2}}{2(1+r(\k\q,\l\boldsymbol{\delta}))}\right) \nonumber\\
&\leq& C (\prod_{i=1}^{d} T_{i}^{-1}\delta_{i}^{1/2})\left(\log \prod_{i=1}^{d} T_{i}\right)^{-\frac{1}{2}-\sum_{i=1}^{d} \frac{1}{2\alpha_{i}}}\nonumber\\
&\times&\sum_{\mbox{\tiny$\begin{array}{c} \k\q, \l\boldsymbol{\delta} \in \O_{\i}, \k\q\neq\l\boldsymbol{\delta}, \i=\mathbf{1,2,\cdots,n} \\ \max\{|l_{j}\delta_{j}-k_{j}q_{j}|,j=1,2,\cdots,d\}\leq \epsilon\end{array}$}} \sqrt{\sum_{j=1}^{d}|l_{j}\delta_{j}-k_{j}q_{j}|^{\alpha_{j}}}  \exp\left(-\frac{u_{\T}^{2}}{8}\left(\sum_{j=1}^{d}|l_{j}\delta_{j}-k_{j}q_{j}|^{\alpha_{j}}\right)\right).
\end{eqnarray}
Recall that $ q_{i}=\gamma_{i}u^{-2/\alpha_{i}}$, and $ \Re(\delta_{i}) $, $ i=1,2,\cdots,d $ are sparse grids, then after some calculation, \eqref{eq6.13} can be bounded by
\begin{eqnarray}
\label{eq6.14}
&& C (\prod_{i=1}^{d} T_{i}^{-1}\delta_{i}^{1/2})\left(\log \prod_{i=1}^{d} T_{i}\right)^{-\frac{1}{2}-\sum_{i=1}^{d} \frac{1}{2\alpha_{i}}}\left(\prod_{i=1}^{d} T_{i} \delta_{i}^{-1}\right)\nonumber\\
&\times&\sum_{j=1}^{d}  \sum_{ 0\leq k_{j}q_{j} \leq \epsilon} \exp\left(-\frac{1}{4}\left[\sum_{i=1}^{d} (k_{i}q_{i})^{\alpha_{i}}\right]  \left(\log \prod_{i=1}^{d} T_{i}\right) \right) \nonumber\\
&\leq& Cu^{-1}\left(\prod_{i=1}^{d}\delta_{i}^{-1/2}u^{-1/\alpha_{i}}\right) \sum_{j=1}^{d} \sum_{k_{j}=0}^{\infty} e^{-\frac{1}{4}(k_{j}\gamma_{j})^{\alpha_{j}}}
\to 0
\end{eqnarray}
as $ u\to\infty $.

Noting that $ u_{\T}\sim \left(2\log\prod_{i=1}^{d} T_{i}\right)^{1/2}=u $, we have
\begin{eqnarray}
\label{eq6.15}
W_{\T,2}&\leq& C \sum_{\mbox{\tiny$\begin{array}{c} \k\q, \l\boldsymbol{\delta} \in \O_{\i}, \k\q\neq\l\boldsymbol{\delta}, \i=\mathbf{1,2,\cdots,n} \\ \max\{|l_{j}\delta_{j}-k_{j}q_{j}|,j=1,2,\cdots,d\}> \epsilon\end{array}$}} | r(\k\q,\l\boldsymbol{\delta})-\varrho(\k\q,\l\boldsymbol{\delta})| \exp\left(-\frac{u_{\T}^{2}}{1+\varpi(\k\q,\l\boldsymbol{\delta})}\right)\nonumber\\
&\leq & C (\prod_{i=1}^{d} T_{i}\delta_{i}^{-1}) \sum_{j=1}^{d} \sum_{\mbox{\tiny$\begin{array}{c} 0\leq k_{j}q_{j}\leq T_{j}^{a} \\ \max\{k_{j}q_{j},j=1,2,\cdots,d\}> \epsilon\end{array}$}} \exp\left(-\frac{u_{\T}^{2}}{1+\theta(\epsilon)}\right)   \nonumber\\
&\leq & C \left(\prod_{i=1}^{d} \delta_{i}^{-1}u^{-2/\alpha_{i}}\right)\exp\left[(a-\frac{1-\theta(\epsilon)}{1+\theta(\epsilon)})u^{2}/2+{\sum_{i=1}^{d} \frac{4}{\alpha_{i}}}\log u\right]\to 0
\end{eqnarray}
as $ u\to\infty $.

Next, we consider the case that $ \k\q\in \O_{\i} $, $ \l\boldsymbol{\delta}\in \O_{\j} $, $ \i\neq \j $. Note that the distance between the points in any two rectangles $ \O_{\i} $ and $ \O_{\j} $ is larger than $ \min\{ T_{i}^{b}, i=1,2,\cdots,d \} $ and $ \varrho(\k\q,\l\boldsymbol{\delta})=\rho(\T)$. Then, $ F_{3} $ is at most
\begin{eqnarray}
\label{eq6.16}
C \sum_{\mbox{\tiny$\begin{array}{c} \k\q\in \O_{\i},\l\boldsymbol{\delta}\in \O_{\j} \\ \mathbf{1\leq i\neq j\leq n}\end{array}$}} | r(\k\q,\l\boldsymbol{\delta})-\rho(\T)| \exp\left(-\frac{u_{\T}^{2}}{1+\varpi(\k\q,\l\boldsymbol{\delta})}\right).
\end{eqnarray}
Split \eqref{eq6.16} into two parts, the first for $ \min\{|k_{j}q_{j}-l_{j}\delta_{j}|,j=1,2,\cdots,d\}> 0 $, and the second for $ \min\{|k_{j}q_{j}-l_{j}\delta_{j}|,j=1,2,\cdots,d\}= 0 $ and denote them  as $ H_{\T,1} $, $ H_{\T,2} $, respectively. Then,
\[ H_{\T,1}=C\sum_{\mbox{\tiny$\begin{array}{c} \k\q\in \O_{\i},\l\boldsymbol{\delta}\in \O_{\j},\mathbf{1\leq i\neq j\leq n} \\ \prod_{m=1}^{d}|k_{m}q_{m}-l_{m}\delta_{m}|\leq {\prod_{i=1}^{d} T_{i}}^{\beta}\end{array}$}}+C\sum_{\mbox{\tiny$\begin{array}{c} \k\q\in \O_{\i},\l\boldsymbol{\delta}\in \O_{\j},\mathbf{1\leq i\neq j\leq n} \\ \prod_{m=1}^{d}|k_{m}q_{m}-l_{m}q_{m}|> \prod_{i=1}^{d} T_{i}^{\beta}\end{array}$}} :=H_{\T,11}+H_{\T,12} .\]
\pzx{By the arguments as used in \eqref{eq6.15}}, we have
\begin{eqnarray}
\label{eq6.17}
H_{\T,11}&\leq & C (\prod_{i=1}^{d} T_{i}\delta_{i}^{-1}) \sum_{j=1}^{d} \sum_{\mbox{\tiny$\begin{array}{c} 0\leq k_{j}q_{j}\leq T_{j} \\ \prod_{m=1}^{d} k_{m}q_{m}\leq {\prod_{i=1}^{d} T_{i}}^{\beta} \end{array}$}} \exp\left(-\frac{u_{\T}^{2}}{1+\theta(\epsilon)}\right)   \nonumber\\
&\leq &  C \left(\prod_{i=1}^{d} \delta_{i}^{-1}u^{-2/\alpha_{i}}\right)\exp\left[(\beta-\frac{1-\theta(\epsilon)}{1+\theta(\epsilon)})u^{2}/2+{\sum_{i=1}^{d} \frac{4}{\alpha_{i}}}\log u\right]\nonumber\\
&\to& 0
\end{eqnarray}
as $u\to\infty  $.

For $ H_{\T,12} $, by \pzx{the same arguments as used in \eqref{eq6.9}}, we have
\begin{eqnarray*}
&&\frac{(\prod_{i=1}^{d} T_{i}^{2}q_{i}^{-1} \delta_{i}^{-1})}{\log\left(\prod_{i=1}^{d} T_{i}\right)}\exp\left(-\frac{u_{\T}^{2}}{1+\theta_{1}(\T^{\beta})}\right)\nonumber\\
&\leq& \frac{(\prod_{i=1}^{d}T_{i}^{2}q_{i}^{-1} \delta_{i}^{-1}) }{\log\left(\prod_{i=1}^{d} T_{i}\right)}\exp\left(-\frac{u_{\T}^{2}}{1+C/\log\left(\prod_{i=1}^{d} T_{i}^{\beta}\right)}\right)\nonumber\\
&=&O(1)\left(\prod_{i=1}^{d} T_{i}\right)^{ 2C/\left(C+\log\left(\prod_{i=1}^{d} T_{i}^{\beta}\right)\right)} \times\left(\prod_{i=1}^{d}\delta_{i} \right)^{(-C)/\left(C+\log\left(\prod_{i=1}^{d} T_{i}^{\beta}\right)\right)} \nonumber \\
&&\times \left( \log \prod_{i=1}^{d} T_{i} \right)^{\left((\sum_{i=1}^{d} \frac{1}{\alpha_{i}}-1)C\right)/\left(C+\log\left(\prod_{i=1}^{d} T_{i}^{\beta}\right)\right)}\nonumber\\
&=&O(1)\exp\left[\frac{C}{C/u^{2}+\beta/2}+2C\left(\sum_{i=1}^{d}2/\alpha_{i}-1\right)\frac{\log u}{C+\beta u^{2}/2}\right]\left(\prod_{i=1}^{d}\delta_{i} u^{\frac{2}{\alpha_{i}}} \right)^{(-C)/\left(C+\beta u^{2}/2\right)}\\
&=&O(1)
\end{eqnarray*}
as $u\to\infty$. Then,
\begin{eqnarray}
\label{eq6.18}
H_{\T,12}&\leq & C\frac{(\prod_{i=1}^{d} T_{i}^{2}q_{i}^{-1} \delta_{i}^{-1})}{\log\left(\prod_{i=1}^{d} T_{i}\right)}\exp\left(-\frac{u_{\T}^{2}}{1+\theta_{1}(\T^{\beta})}\right) \frac{\log\left(\prod_{i=1}^{d} T_{i}\right)}{(\prod_{i=1}^{d} T_{i}^{2}q_{i}^{-1} \delta_{i}^{-1})} \nonumber\\ &\times&\sum_{\mbox{\tiny$\begin{array}{c} \k\q\in \O_{\i},\l\boldsymbol{\delta}\in \O_{\j},\mathbf{1\leq i\neq j\leq n} \\ \prod_{m=1}^{d}|k_{m}q_{m}-l_{m}q_{m}|> \prod_{i=1}^{d} T_{i}^{\beta}\end{array}$}} |r(\k\q-\l\boldsymbol{\delta})-\rho(\T) | \nonumber\\
&\leq& C \left(\prod_{i=1}^{d} \delta_{i} q_{i} T_{i}^{-2}\right)\sum_{\mbox{\tiny$\begin{array}{c} \k\q\in \O_{\i},\l\boldsymbol{\delta}\in \O_{\j},\mathbf{1\leq i\neq j\leq n} \\ \prod_{m=1}^{d}|k_{m}q_{m}-l_{m}q_{m}|> \prod_{i=1}^{d} T_{i}^{\beta}\end{array}$}} \left|r(\k\q-\l\boldsymbol{\delta})\log\left(\prod_{m=1}^{d} k_{m}q_{m}-l_{m}\delta_{m}\right)-r \right| \nonumber\\
&+& C \left(\prod_{i=1}^{d} \delta_{i} q_{i} T_{i}^{-2}\right)\sum_{\mbox{\tiny$\begin{array}{c} \k\q\in \O_{\i},\l\boldsymbol{\delta}\in \O_{\j},\mathbf{1\leq i\neq j\leq n} \\ \prod_{m=1}^{d}|k_{m}q_{m}-l_{m}\delta_{m}|> \prod_{i=1}^{d} T_{i}^{\beta}\end{array}$}} \left| 1-\frac{\log\prod_{i=1}^{d} T_{i}}{\log\prod_{j=1}^{d} (k_{m}q_{m}-l_{m}\delta_{m})} \right|.
\end{eqnarray}
By condition $ \mathbf{A3} $, the first term on the right hand side of \eqref{eq6.18} tends to $ 0 $ as $ \T\to\infty $, and the second term also tends to $ 0 $ by the same arguments of \eqref{eq6.11}, then $ H_{\T,12}\to 0 $ as $ \T\to\infty $. Combining \eqref{eq6.17}, we can get that $ H_{\T,1}\to 0 $ as $ \T\to\infty $.

Now, we consider $ H_{\T,2} $, we only prove the case that $ k_{1}q_{1}= l_{1}\delta_{1} $, and $ \min\{|k_{j}q_{j}-l_{j}\delta_{j}|,j=2,\cdots,d\}> 0  $. By the similar way, the rest cases can be proved. If $ \prod_{i=2}^{d} T_{i}\leq \prod_{i=1}^{d} T_{i}^{\beta}=e^{\beta u^{2}/2} $, by \pzx{the same arguments as used in \eqref{eq6.8} that}
\begin{eqnarray}
\label{eq6.19}
H_{\T,2}&\leq & C  \sum_{\mbox{\tiny$\begin{array}{c} \k\q\in \O_{\i}, \l\boldsymbol{\delta}\in\O_{\j}, \mathbf{1\leq i\neq j\leq n} \\ k_{1}q_{1}=l_{1}\delta_{1}\end{array}$}} |r(\k\q-\l\boldsymbol{\delta})-\rho(\T) |\exp\left(-\frac{u_{\T}^{2}(x_{2})}{1+\theta(\epsilon)}\right)   \nonumber\\
&\leq &  C \prod_{i=1}^{d} \delta_{i}^{-1}u^{-2/\alpha_{i}}\exp\left[(\beta-\frac{1-\theta(\epsilon)}{1+\theta(\epsilon)})u^{2}/2+{\sum_{i=1}^{d} \frac{4}{\alpha_{i}}}\log u\right]\to 0
\end{eqnarray}
as $ u\to\infty $. If $ \prod_{i=2}^{d} T_{i}> \prod_{i=1}^{d} T_{i}^{\beta}=e^{\beta u^{2}/2} $, we have
\[ H_{\T,2}=C\sum_{\mbox{\tiny$\begin{array}{c} \k\q\in \O_{\i},\l\boldsymbol{\delta}\in \O_{\j},\mathbf{1\leq i\neq j\leq n} \\ \prod_{m=2}^{d}|k_{m}q_{m}-l_{m}q_{m}|\leq {\prod_{i=1}^{d} T_{i}}^{\beta}\\ k_{1}q_{1}=l_{1}\delta_{1}\end{array}$}}+C\sum_{\mbox{\tiny$\begin{array}{c} \k\q\in \O_{\i},\l\boldsymbol{\delta}\in \O_{\j},\mathbf{1\leq i\neq j\leq n} \\ \prod_{m=2}^{d}|k_{m}q_{m}-l_{m}q_{m}|> \prod_{i=1}^{d} T_{i}^{\beta}\\ k_{1}q_{1}=l_{1}\delta_{1}\end{array}$}} :=H_{\T,21}+H_{\T,22} .\]
For $ H_{\T,21} $, it follows from \eqref{eq6.19} that
\begin{eqnarray*}
H_{\T,21}&\leq & C \sum_{\mbox{\tiny$\begin{array}{c} \k\q\in \O_{\i},\l\boldsymbol{\delta}\in \O_{\j},\mathbf{1\leq i\neq j\leq n} \\ \prod_{m=2}^{d}|k_{m}q_{m}-l_{m}q_{m}|\leq {\prod_{i=1}^{d} T_{i}}^{\beta}\\ k_{1}q_{1}=l_{1}\delta_{1}\end{array}$}}  |r(\k\q-\l\boldsymbol{\delta})-\rho(\T) |\exp\left(-\frac{u_{\T}^{2}(x_{2})}{1+\theta(\epsilon)}\right)   \nonumber\\
&\to& 0\;\ as\;\ \T\to\infty.
\end{eqnarray*}
For $ H_{\T,22} $, \pzx{by the same arguments as used in \eqref{eq6.9}}, we have
\[ \frac{(\prod_{i=1}^{d} T_{i}^{2}q_{i}^{-1}\delta_{i}^{-1})}{\log\left(\prod_{i=1}^{d} T_{i}\right)}\exp\left(-\frac{u_{\T}^{2}}{1+\theta_{2}(\T^{\beta})}\right)= O(1),\]
hence,
\begin{eqnarray*}
H_{\T,22}&\leq & C \sum_{\mbox{\tiny$\begin{array}{c} \k\q\in \O_{\i},\l\boldsymbol{\delta}\in \O_{\j},\mathbf{1\leq i\neq j\leq n} \\ \prod_{m=2}^{d}|k_{m}q_{m}-l_{m}q_{m}|> \prod_{i=1}^{d} T_{i}^{\beta}\\ k_{1}q_{1}=l_{1}\delta_{1}\end{array}$}} \left| r(0,k_{2}q_{2},\cdots,k_{d}q_{d})-\rho(\T)  \right| \exp\left(-\frac{u_{\T}^{2}}{1+\theta_{2}(\T^{\beta})}\right)\\
&= & C \frac{(\prod_{i=1}^{d} T_{i}^{2}q_{i}^{-1}\delta_{i}^{-1})}{\log\left(\prod_{i=1}^{d} T_{i}\right)}\exp\left(-\frac{u_{\T}^{2}}{1+\theta_{2}(\T^{\beta})}\right)
\times\frac{\log(\prod_{i=1}^{d} T_{i})}{(\prod_{i=1}^{d} T_{i}^{2}q_{i}^{-1}\delta_{i}^{-1})} \\ &\times&\sum_{\mbox{\tiny$\begin{array}{c} \k\q\in \O_{\i},\l\boldsymbol{\delta}\in \O_{\j},\mathbf{1\leq i\neq j\leq n} \\ \prod_{m=2}^{d}|k_{m}q_{m}-l_{m}q_{m}|> \prod_{i=1}^{d} T_{i}^{\beta}\\ k_{1}q_{1}=l_{1}\delta_{1}\end{array}$}} \left| r(0,k_{2}q_{2},\cdots,k_{d}q_{d})-\rho(\T)  \right|\nonumber\\
&\leq& C\frac{\log(\prod_{i=1}^{d} T_{i})}{(\prod_{i=1}^{d} T_{i}^{2}q_{i}^{-1}\delta_{i}^{-1})} \sum_{\mbox{\tiny$\begin{array}{c} \k\q\in \O_{\i},\l\boldsymbol{\delta}\in \O_{\j},\mathbf{1\leq i\neq j\leq n} \\ \prod_{m=2}^{d}|k_{m}q_{m}-l_{m}q_{m}|> \prod_{i=1}^{d} T_{i}^{\beta}\\ k_{1}q_{1}=l_{1}\delta_{1}\end{array}$}} \left( |r(0,k_{2}q_{2},\cdots,k_{d}q_{d})|+|\rho(\T)|  \right)\\
&\leq & C (\prod_{i=1}^{d} T_{i}^{-2}q_{i}\delta_{i}) \left(\prod_{i=1}^{d} T_{i}\delta_{i}^{-1}\right)\left(\prod_{i=2}^{d} T_{i}q_{i}^{-1}\right)\nonumber\\
&\leq& C\frac{q_{1}}{T_{1}},
\end{eqnarray*}
which implies $ H_{\T,22}\to 0 $ as $ \T\to\infty $. Furthermore, $ H_{\T,2}\to 0 $, as $ \T\to\infty $. Combining \eqref{eq6.14}, \eqref{eq6.15}, \eqref{eq6.17} and \eqref{eq6.18}, we can get $ F_{3}\to 0 $ as $ \T\to\infty $. The proof  is complete.
\qed

\begin{lemma}
\label{lem7}
Suppose $ \Re(\delta_{i}), i=1,2,\cdots,d $ are all sparse grids or all Pickands grids and  ${\u}_{\T}$, $\v_{\T}$ are given by \eqref{eq4.2}.  Then for the grids $ \Re(q_{i})$ with $ q_{i}=\gamma_{i}u^{-2/\alpha_{i}} $  and $ \gamma_{i}>0 $, $i=1,2,\cdots,d $,  we have
\begin{eqnarray}
\label{eq7.1}
&&\Big|\P\Big( v_{\T}(x_{1}) < \min_{\t\in \O\cap\prod_{i=1}^{d} \Re(q_{i})}\xi_{\T}(\t) \leq \max_{\t\in \O\cap \prod_{i=1}^{d} \Re(q_{i})} \xi_{\T}(\t)\leq u_{\T}(x_{2}), \nonumber\\
&& \qquad v_{\T}^{\boldsymbol{\delta}}(y_{1}) < \min_{\t\in \O \cap\prod_{i=1}^{d} \Re(\delta_{i})} \xi_{\T}(\t) \leq \max_{\t\in \O\cap\prod_{i=1}^{d} \Re(\delta_{i})} \xi_{\T}(\t)\leq  u_{\T}^{\boldsymbol{\delta}}(y_{2})\Big) \nonumber\\
&&\qquad\qquad-\int_{-\infty}^{+\infty} \prod_{j=1}^{d} \prod_{i_{j}=1}^{n_{j}}\P\Big( v_{\T}^{'}(x_{1})<\min_{\t\in\O_{\i}} \eta(\t)\leq \max_{\t\in\O_{\i}} \eta(\t)\leq u_{\T}^{'}(x_{2}),\nonumber\\
&& \qquad\qquad\qquad v_{\T}^{*}(y_{1})<\min_{\t\in\O_{\i}\cap\prod_{i=1}^{d} \Re(\delta_{i})} \eta(\t)\leq \max_{\t\in\O_{\i}\cap\prod_{i=1}^{d} \Re(\delta_{i})} \eta(\t)\leq u_{\T}^{*}(y_{2}) \Big) \phi(z) dz \Big| \nonumber\\
&&\to 0
\end{eqnarray}
as $\T\to\infty$ and all $ \gamma_{i}\downarrow 0 $,  where
\begin{eqnarray}
\label{eq7.2}
\begin{cases}
v_{\T}^{'}(x_{1})=\frac{-b_{\T}+x_{1}/a_{\T}-\rho^{1/2}(\T)z}{(1-\rho(\T))^{1/2}}=\frac{x_{1}-r-\sqrt{2r}z}{a_{\T}}-b_{\T}+o(a_{\T}^{-1});\\
u_{\T}^{'}(x_{2})=\frac{b_{\T}+x_{2}/a_{\T}-\rho^{1/2}(\T)z}{(1-\rho(\T))^{1/2}}=\frac{x_{2}+r-\sqrt{2r}z}{a_{\T}}+b_{\T}+o(a_{\T}^{-1});\\
v_{\T}^{*}(y_{1})=\frac{-b_{\T}^{*}+y_{1}/a_{\T}-\rho^{1/2}(\T)z}{(1-\rho(\T))^{1/2}}=\frac{y_{1}-r-\sqrt{2r}z}{a_{\T}}-b_{\T}^{*}+o(a_{\T}^{-1});\\
u_{\T}^{*}(y_{2})=\frac{b_{\T}^{*}+y_{2}/a_{\T}-\rho^{1/2}(\T)z}{(1-\rho(\T))^{1/2}}=\frac{y_{2}+r-\sqrt{2r}z}{a_{\T}}+b_{\T}^{*}+o(a_{\T}^{-1}).
\end{cases}
\end{eqnarray}
\end{lemma}

\noindent {\bf Proof.} By the definition of $ \{\xi_{\T}(\t),  \t\geq \mathbf{0} \} $ and $\{ \eta(\t),\t\geq \mathbf{0}\} $, we have
\begin{eqnarray}
\label{eq7.3}
&&\P\Big( v_{\T}(x_{1}) < \min_{\t\in \O\cap\prod_{i=1}^{d} \Re(q_{i})}\xi_{\T}(\t) \leq \max_{\t\in \O\cap \prod_{i=1}^{d} \Re(q_{i})} \xi_{\T}(\t)\leq u_{\T}(x_{2}), \nonumber\\
&&\qquad v_{\T}^{\boldsymbol{\delta}}(y_{1}) < \min_{\t\in \O \cap\prod_{i=1}^{d} \Re(\delta_{i})} \xi_{\T}(\t) \leq \max_{\t\in \O\cap\prod_{i=1}^{d} \Re(\delta_{i})} \xi_{\T}(\t)\leq  u_{\T}^{\boldsymbol{\delta}}(y_{2})\Big) \nonumber\\
&=&\int_{-\infty}^{+\infty} \P\Big( v_{\T}^{'}(x_{1}) < \min_{\t\in \O\cap\prod_{i=1}^{d} \Re(q_{i})}\eta(\t) \leq \max_{\t\in \O \cap\prod_{i=1}^{d} \Re(q_{i})} \eta(\t)\leq u_{\T}^{'}(x_{2}), \nonumber\\ && v_{\T}^{*}(y_{1}) < \min_{\t\in \O \cap\prod_{i=1}^{d} \Re(\delta_{i})} \eta(\t) \leq \max_{\t\in \O\cap\prod_{i=1}^{d} \Re(\delta_{i})} \eta(\t)\leq  u_{\T}^{*}(y_{2}) \Big) \phi(z) dz \nonumber\\
&=&\int_{-\infty}^{+\infty} \prod_{j=1}^{d} \prod_{i_{j}=1}^{n_{j}}\P\Big( v_{\T}^{'}(x_{1}) < \min_{\t\in \O_{\i}\cap\prod_{i=1}^{d} \Re(q_{i})}\eta(\t) \leq \max_{\t\in \O_{\i}\cap \prod_{i=1}^{d} \Re(q_{i})} \eta(\t)\leq u_{\T}^{'}(x_{2}), \nonumber\\
 && \qquad v_{\T}^{*}(y_{1}) < \min_{\t\in \O_{\i} \cap\prod_{i=1}^{d} \Re(\delta_{i})} \eta(\t) \leq \max_{\t\in \O_{\i}\cap\prod_{i=1}^{d} \Re(\delta_{i})} \eta(\t)\leq  u_{\T}^{*}(y_{2}) \Big) \phi(z) dz.
\end{eqnarray}
By Lemma \ref{lem3} and the dominated convergence theorem, we have
\begin{eqnarray*}
&&\Big|\int_{-\infty}^{+\infty} \prod_{j=1}^{d} \prod_{i_{j}=1}^{n_{j}}\P\Big( v_{\T}^{'}(x_{1}) < \min_{\t\in \O_{\i}\cap\prod_{i=1}^{d} \Re(q_{i})}\eta(\t) \leq \max_{\t\in \O_{\i}\cap \prod_{i=1}^{d} \Re(q_{i})} \eta(\t)\leq u_{\T}^{'}(x_{2}), \nonumber\\ && v_{\T}^{*}(y_{1}) < \min_{\t\in \O_{\i} \cap\prod_{i=1}^{d} \Re(\delta_{i})} \eta(\t) \leq \max_{\t\in \O_{\i}\cap\prod_{i=1}^{d} \Re(\delta_{i})} \eta(\t)\leq  u_{\T}^{*}(y_{2}) \Big) \phi(z) dz \nonumber\\
&&-\int_{-\infty}^{+\infty} \prod_{j=1}^{d} \prod_{i_{j}=1}^{n_{j}}\P\Big( v_{\T}^{'}(x_{1})<\min_{\t\in\O_{\i}} \eta(\t)\leq \max_{\t\in\O_{\i}} \eta(\t)\leq u_{\T}^{'}(x_{2}),\nonumber\\ && v_{\T}^{*}(y_{1}) < \min_{\t\in \O_{\i} \cap\prod_{i=1}^{d} \Re(\delta_{i})} \eta(\t) \leq \max_{\t\in \O_{\i}\cap\prod_{i=1}^{d} \Re(\delta_{i})} \eta(\t)\leq  u_{\T}^{*}(y_{2}) \Big) \phi(z) dz \Big| \\
&&\to 0
\end{eqnarray*}
as $\T\to\infty$ and all $ \gamma_{i}\downarrow 0 $. Combining \eqref{eq7.3}, we finish the proof.
\qed

\begin{lemma}
\label{lem8}
Let $ X(\t) $ be a centered homogeneous Gaussian field with \pzx{unit variance} and covariance function $ r(\t) $ satisfying $ \mathbf{A1}-\mathbf{A3} $. Then for any sparse grids,
\begin{eqnarray*}
&&\lim_{\T\to\infty}\P\left(\mathbf{v}_{\T}<\mathbf{m}_{\T}\leq \mathbf{M}_{\T}\leq \mathbf{u}_{\T}\right)\\
&=&\int_{-\infty}^{+\infty} \exp\Big( -(e^{-x_{2}-r+\sqrt{2r}z}
+e^{-y_{2}-r+\sqrt{2r}z}+e^{x_{1}+r+\sqrt{2r}z}+e^{y_{1}+r+\sqrt{2r}z})\Big) \phi(z) dz,
\end{eqnarray*}
where $\mathbf{u}_{\T}$, $\mathbf{v}_{\T}$ are given by \eqref{eq4.2}.
\end{lemma}
\noindent {\bf Proof.} Noting that $\{\eta(\t),\t\geq \mathbf{0}\}=^{d} \{-\eta(\t),\t\geq \mathbf{0}\} $, we have
\begin{eqnarray}
\label{eq8.1}
&&1-\P\Big(v_{\T}^{'}(x_{1})<\min_{\t\in\T_{\T^{a}}} \eta(\t)\leq \max_{\t\in\I_{\T^{a}}} \eta(\t)\leq u_{\T}^{'}(x_{2}),\nonumber\\
&&\qquad v_{\T}^{*}(y_{1}) < \min_{\t\in \I_{\T^{a}} \cap\prod_{i=1}^{d} \Re(\delta_{i})} \eta(\t) \leq \max_{\t\in \I_{\T^{a}}\cap\prod_{i=1}^{d} \Re(\delta_{i})} \eta(\t)\leq  u_{\T}^{*}(y_{2}) \Big) \nonumber\\
&=& 1-\P\left(\max_{\t\in\I_{\T^{a}}} \eta(\t)\leq u_{\T}^{'}(x_{2}), \max_{\t\in \I_{\T^{a}}\cap\prod_{i=1}^{d} \Re(\delta_{i})} \eta(\t)\leq  u_{\T}^{*}(y_{2}) \right)\nonumber\\
&&+1-\P\left(\max_{\t\in\I_{\T^{a}}} -\eta(\t)<-v_{\T}^{'}(x_{1}),\max_{\t\in \I_{\T^{a}} \cap\prod_{i=1}^{d} \Re(\delta_{i})} -\eta(\t)< -v_{\T}^{*}(y_{1}) \right)\nonumber\\
&& -\P\left(\max_{\t\in\I_{\T^{a}}} \eta(\t)> u_{\T}^{'}(x_{2}),\max_{\t\in \I_{\T^{a}} \cap\prod_{i=1}^{d} \Re(\delta_{i})} -\eta(\t)\geq -v_{\T}^{*}(y_{1}) \right)\nonumber\\
&&-\P\left(\max_{\t\in\I_{\T^{a}}} -\eta(\t)\geq-v_{\T}^{'}(x_{1}), \max_{\t\in \I_{\T^{a}}\cap\prod_{i=1}^{d} \Re(\delta_{i})} \eta(\t)>  u_{\T}^{*}(y_{2})  \right).
\end{eqnarray}
Since
\begin{eqnarray*}
&&\P\left(\max_{\t\in\I_{\T^{a}}} \eta(\t)> u_{\T}^{'}(x_{2}),\max_{\t\in \I_{\T^{a}} \cap\prod_{i=1}^{d} \Re(\delta_{i})} -\eta(\t)\geq -v_{\T}^{*}(y_{1}) \right)\\
&&\leq \P\left(\max_{(\t,\s)\in\I_{\T^{a}}\times \I_{\T^{a}}} \eta(\t)-\eta(\s)> u_{\T}^{'}(x_{2})-v_{\T}^{*}(y_{1}) \right)\\
&&\leq 2\P\left( \max_{\t\in\I_{\T^{a}}} \eta(\t)>\frac{u_{\T}^{'}(x_{2})-v_{\T}^{*}(y_{1})}{2} \right),
\end{eqnarray*}
by Borell theorem and \eqref{eq7.2} that for some constant $ c $,  with \pzx{$ n_{i}=\lfloor T_{i}/(T_{i}^{a}+T_{i}^{b})\rfloor$} defined as before we have{\small
\begin{eqnarray}
\label{eq8.2}
&& \prod_{i=1}^{d} n_{i} \P\left(\max_{\t\in\I_{\T^{a}}} \eta(\t)> u_{\T}^{'}(x_{2}),\max_{\t\in \I_{\T^{a}} \cap\prod_{i=1}^{d} \Re(\delta_{i})} -\eta(\t)\geq -v_{\T}^{*}(y_{1}) \right)\nonumber\\
& \leq &C \prod_{i=1}^{d} n_{i} \Psi\left(\frac{u_{\T}^{'}(x_{2})-v_{\T}^{*}(y_{1})}{2} -c\right)\nonumber\\
& \leq &C \prod_{i=1}^{d} n_{i} \left(2\log \prod_{i=1}^{d} T_{i}\right)^{-\frac{1}{2}}\Big(1+\frac{-2c\sqrt{2\log \prod_{i=1}^{d} T_{i}}+x_{2}-y_{1}+2r}{4\log\left(\prod_{i=1}^{d} T_{i}\right)}+\nonumber\\
&&\frac{\log\left[(2\pi)^{-1/2}\left(\prod_{i=1}^{d} H_{\alpha_{i}}\right)\left(2\log \prod_{i=1}^{d} T_{i}\right)^{-1/2+\sum_{i=1}^{d}\frac{1}{\alpha_{i}}}\right]+\log\left[(2\pi)^{-1/2}\left(\prod_{i=1}^{d} \delta_{i}^{-1}\right)\left(2\log \prod_{i=1}^{d} T_{i}\right)^{-1/2}\right]}{4\log\left(\prod_{i=1}^{d} T_{i}\right)}\Big)^{-1}\nonumber\\
&&\times \exp\Big( -\log\left(\prod_{i=1}^{d} T_{i}\right)-\frac{-2c\sqrt{2\log \prod_{i=1}^{d} T_{i}}+x_{2}-y_{1}+2r}{2} -\nonumber\\
&&\frac{\log\left[(2\pi)^{-1/2}\left(\prod_{i=1}^{d} H_{\alpha_{i}}\right)\left(2\log \prod_{i=1}^{d} T_{i}\right)^{-1/2+\sum_{i=1}^{d}\frac{1}{\alpha_{i}}}\right]+\log\left[(2\pi)^{-1/2}\left(\prod_{i=1}^{d} \delta_{i}^{-1}\right)\left(2\log \prod_{i=1}^{d} T_{i}\right)^{-1/2}\right]}{2} \Big)\nonumber\\
& \leq &C \exp\Big(-\sum_{i=1}^{d} \log(T_{i}^{a}+T_{i}^{b})-\sum_{i=1}^{d}\frac{1}{2\alpha_{i}}\log(2\log \prod_{i=1}^{d} T_{i})+c\sqrt{2\log (\prod_{i=1}^{d} T_{i})}\Big)
\to 0
\end{eqnarray}
}
as $\T\to\infty.$ Similarly,
\begin{eqnarray}
\label{eq8.3}
\prod_{i=1}^{d} n_{i} \P\left(\max_{\t\in\T_{\T^{a}}} -\eta(\t)\geq-v_{\T}^{'}(x_{1}), \max_{\t\in \I_{\T^{a}}\cap\prod_{i=1}^{d} \Re(\delta_{i})} \eta(\t)>  u_{\T}^{*}(y_{2})\right) \to 0
\end{eqnarray}
as $ \T\to\infty $. Using Lemma \ref{lem2} and \eqref{eq7.1}, we have
\begin{eqnarray}
\label{eq8.4}
&&\prod_{i=1}^{d} n_{i} \left(1-\P\left(\max_{\t\in\I_{\T^{a}}} \eta(\t)\leq u_{\T}^{'}(x_{2}), \max_{\t\in \I_{\T^{a}}\cap\prod_{i=1}^{d} \Re(\delta_{i})} \eta(\t)\leq  u_{\T}^{*}(y_{2}) \right)\right) \nonumber\\
&\to& e^{-x_{2}-r+\sqrt{2r}z}+e^{-y_{2}-r+\sqrt{2r}z}
\end{eqnarray}
as $ \T\to\infty $, and
\begin{eqnarray}
\label{eq8.5}
&&\prod_{i=1}^{d} n_{i} \left(1-\P\left(\max_{\t\in\T_{\T^{a}}} -\eta(\t)<-v_{\T}^{'}(x_{1}),\max_{\t\in \I_{\T^{a}} \cap\prod_{i=1}^{d} \Re(\delta_{i})} -\eta(\t)< -v_{\T}^{*}(y_{1}) \right)\right)\nonumber\\
&\to& e^{x_{1}+r+\sqrt{2r}z}+e^{y_{1}+r+\sqrt{2r}z}
\end{eqnarray}
as $ \T\to\infty $. Since
\begin{eqnarray}
\label{eq8.6}
&&\prod_{j=1}^{d} \prod_{i_{j}=1}^{n_{j}}\P\Big( v_{\T}^{'}(x_{1})<\min_{\t\in\O_{\i}} \eta(\t)\leq \max_{\t\in\O_{\i}} \eta(\t)\leq u_{\T}^{'}(x_{2}),\nonumber\\
&&\qquad v_{\T}^{*}(y_{1}) < \min_{\t\in \O_{\i} \cap\prod_{i=1}^{d} \Re(\delta_{i})} \eta(\t) \leq \max_{\t\in \O_{\i}\cap\prod_{i=1}^{d} \Re(\delta_{i})} \eta(\t)\leq  u_{\T}^{*}(y_{2}) \Big)\nonumber\\
&=&\Big( \P\Big(v_{\T}^{'}(x_{1})<\min_{\t\in\I_{\T^{a}}} \eta(\t)\leq \max_{\t\in\I_{\T^{a}}} \eta(\t)\leq u_{\T}^{'}(x_{2}),\nonumber\\
&& v_{\T}^{*}(y_{1}) < \min_{\t\in \I_{\T^{a}} \cap\prod_{i=1}^{d} \Re(\delta_{i})} \eta(\t) \leq \max_{\t\in \I_{\T^{a}}\cap\prod_{i=1}^{d} \Re(\delta_{i})} \eta(\t)\leq  u_{\T}^{*}(y_{2}) \Big)\Big)^{\prod_{i=1}^{d} n_{i}}\nonumber\\
&&=\exp\Big(\prod_{i=1}^{d} n_{i}\log\Big(\P\Big(v_{\T}^{'}(x_{1})<\min_{\t\in\I_{\T^{a}}} \eta(\t)\leq \max_{\t\in\I_{\T^{a}}} \eta(\t)\leq u_{\T}^{'}(x_{2}),\nonumber\\
&&\qquad v_{\T}^{*}(y_{1}) < \min_{\t\in \I_{\T^{a}} \cap\prod_{i=1}^{d} \Re(\delta_{i})} \eta(\t) \leq \max_{\t\in \I_{\T^{a}}\cap\prod_{i=1}^{d} \Re(\delta_{i})} \eta(\t)\leq  u_{\T}^{*}(y_{2}) \Big)\Big)\Big)\nonumber\\
&=&\exp\left(-\prod_{i=1}^{d} n_{i}(1-P_{\T})+R_{\T}\right),
\end{eqnarray}
where
\begin{eqnarray*}
&&P_{\T}=\P\Big(v_{\T}^{'}(x_{1})<\min_{\t\in\I_{\T^{a}}} \eta(\t)\leq \max_{\t\in\I_{\T^{a}}} \eta(\t)\leq u_{\T}^{'}(x_{2}),\nonumber\\
&&\qquad v_{\T}^{*}(y_{1}) < \min_{\t\in \I_{\T^{a}} \cap\prod_{i=1}^{d} \Re(\delta_{i})} \eta(\t) \leq \max_{\t\in \I_{\T^{a}}\cap\prod_{i=1}^{d} \Re(\delta_{i})} \eta(\t)\leq  u_{\T}^{*}(y_{2}) \Big)\\
&&\to 1
\end{eqnarray*}
as $ \T\to\infty $, and
\begin{eqnarray*}
R_{\T}=o\left(\prod_{i=1}^{d} n_{i}(1-P_{\T})\right)
\end{eqnarray*}
as $ \T\to\infty $. Then, combining \eqref{eq8.1}-\eqref{eq8.6} and Lemmas \ref{lem4}-\ref{lem7}, the assertion of this lemma follows.

\begin{lemma}
\label{lem9}
Let $ X(\t) $ be a centered homogeneous Gaussian field with \pzx{unit variance} and covariance function $ r(\t) $ satisfying $ \mathbf{A1}-\mathbf{A3} $ and the $\mathbf{u}_{\T}$, $\mathbf{v}_{\T}$ be given by \eqref{eq4.2}. Then \pzx{for} $H_{\a,\boldsymbol{\alpha}}^{x,y}$ given \pzx{in} Theorem \ref{th2} and any Pickands grids $\Re(a_{i}(2\log \prod_{i=1}^{d} T_{i})^{-1/\alpha_{i}})$ with $a_{i}>0$,  $ i=1,2,\cdots,d $, we have
\begin{eqnarray*}
&&\lim_{\T\to\infty}\P\left(\mathbf{v}_{\T}<\mathbf{m}_{\T}\leq \mathbf{M}_{\T}\leq \mathbf{u}_{\T}\right)\\
&&=\int_{-\infty}^{+\infty}\exp\Big( -\Big(e^{-x_{2}-r+\sqrt{2r}z}+e^{-y_{2}-r+\sqrt{2r}z}
-H_{\a,\boldsymbol{\alpha}}^{x_{2}+\log( \prod_{i=1}^{d}H_{\alpha_{i}}),y_{2}+\log(\prod_{i=1}^{d} H_{a_{i},\alpha_{i}})}e^{-r+\sqrt{2r}z}\\&&+e^{x_{1}+r+\sqrt{2r}z}
+e^{y_{1}+r+\sqrt{2r}z}
-H_{\a,\boldsymbol{\alpha}}^{-x_{1}+\log( \prod_{i=1}^{d}H_{\alpha_{i}}),-y_{1}+\log(\prod_{i=1}^{d} H_{a_{i},\alpha_{i}})}e^{r+\sqrt{2r}z}\Big) \Big)\phi(z)dz.
\end{eqnarray*}
\end{lemma}

\noindent {\bf Proof.} Suppose $ \Re(\delta_{i})$ are all Pickands grids and $S_{i}$ satisfies the conditions in Lemma \ref{lem3}. Denote
\[ P_{\S}^{'}(u,x)=\P\left( \max_{\t\in\I_{\S}\cap\prod_{i=1}^{d} \Re(\delta_{i})} X(\t)> u, \max_{\t\in\I_{\S}} X(\t)> u+\frac{x}{u}\right),\]
then by using Lemma 6.1 of Piterbarg (1996), we have
\begin{eqnarray}
\label{eq2}
P_{\S}^{'}(u,x)\sim \left(\prod_{i=1}^{d} S_{i}u^{2/\alpha_{i}}\right) H_{\a,\boldsymbol{\alpha}}^{0,x} \Psi(u)
\end{eqnarray}
as $ u\to\infty $. Note that for Pickands grids,
\begin{eqnarray*}
u_{\T}^{'}(x_{2})
&=&  u_{\T}^{*}(y_{2})+\frac{\log( \prod_{i=1}^{d}H_{\alpha_{i}})-\log(\prod_{i=1}^{d} H_{a_{i},\alpha_{i}})+x_{2}-y_{2}}{u_{\T}^{*}(y_{2})}
+O\left(\frac{\log u}{u^{3}}\right)
\end{eqnarray*}
and
\begin{eqnarray*}
v_{\T}^{'}(x_{1})
&=&  v_{\T}^{*}(y_{1})+\frac{\log( \prod_{i=1}^{d}H_{\alpha_{i}})-\log(\prod_{i=1}^{d} H_{a_{i},\alpha_{i}})-x_{1}+y_{1}}{v_{\T}^{*}(y_{1})}
+O\left(\frac{\log u}{u^{3}}\right).
\end{eqnarray*}

Hence, from \eqref{eq2}, we have
\begin{eqnarray}
\label{eq9.2}
&&\prod_{i=1}^{d} n_{i} \left(\P\left(\max_{\t\in\I_{\T^{a}}} \eta(\t)> u_{\T}^{'}(x_{2}), \max_{\t\in \I_{\T^{a}}\cap\prod_{i=1}^{d} \Re(\delta_{i})} \eta(\t)>  u_{\T}^{*}(y_{2}) \right)\right) \nonumber\\
&&\sim \left(\prod_{i=1}^{d} n_{i} T_{i}^{a}\right) H_{\a,\boldsymbol{\alpha}}^{0,Z_{x_{2},y_{2}}}(u_{\T}^{*}(y_{2}))^{\sum_{i=1}^{d} \frac{2}{\alpha_{i}}}\Psi(u_{\T}^{*}(y_{2}))\nonumber\\
&&\sim H_{\a,\boldsymbol{\alpha}}^{0,Z_{x_{2},y_{2}}} \left(\prod_{i=1}^{d}H_{a_{i},\alpha_{i}} \right)^{-1}e^{-y_{2}-r+\sqrt{2r}z},
\end{eqnarray}
where $Z_{x_{2},y_{2}}=\log( \prod_{i=1}^{d}H_{\alpha_{i}})-\log(\prod_{i=1}^{d} H_{a_{i},\alpha_{i}})+x_{2}-y_{2} $. By the definition of $ H_{\a,\boldsymbol{\alpha}}^{x,y} $, we get $ H_{\a,\boldsymbol{\alpha}}^{0,Z_{x_{2},y_{2}}}\left(\prod_{i=1}^{d}H_{a_{i},\alpha_{i}} \right)^{-1}e^{-y_{2}}=H_{\a,\boldsymbol{\alpha}}^{x_{2}+\log( \prod_{i=1}^{d}H_{\alpha_{i}}),y_{2}+\log(\prod_{i=1}^{d} H_{a_{i},\alpha_{i}})} $.

Furthermore, from \eqref{eq3.2} and \eqref{eq9.2}, we have
\begin{eqnarray}
\label{eq9.3}
&&\prod_{i=1}^{d} n_{i} \left(1- \P\left(\max_{\t\in\I_{\T^{a}}} \eta(\t)\leq u_{\T}^{'}(x_{2}), \max_{\t\in \I_{\T^{a}}\cap\prod_{i=1}^{d} \Re(\delta_{i})} \eta(\t)\leq  u_{\T}^{*}(y_{2}) \right) \right) \nonumber\\
&&=\prod_{i=1}^{d} n_{i}\Big[\P\left( \max_{\t\in\I_{\T^{a}}} \eta(\t)> u_{\T}^{'}(x_{2})\right)+\P\left( \max_{\t\in\I_{\T^{a}}\cap\prod_{i=1}^{d} \Re(\delta_{i})} \eta(\t)> u_{\T}^{*}(y_{2})\right)\nonumber\\
&&-\P\left( \max_{\t\in\I_{\T^{a}}\cap\prod_{i=1}^{d} \Re(\delta_{i})} \eta(\t)>u_{\T}^{*}(y_{2}) , \max_{\t\in\I_{\T^{a}}} \eta(\t)> u_{\T}^{'}(x_{2})\right)\Big]\nonumber\\
&&\sim e^{-x_{2}-r+\sqrt{2r}z}+e^{-y_{2}-r+\sqrt{2r}z}-H_{\a,\boldsymbol{\alpha}}^{x_{2}+\log( \prod_{i=1}^{d}H_{\alpha_{i}}),y_{2}+\log(\prod_{i=1}^{d} H_{a_{i},\alpha_{i}})}e^{-r+\sqrt{2r}z}.
\end{eqnarray}
Similarly,
\begin{eqnarray}
\label{eq9.4}
&&\prod_{i=1}^{d} n_{i} \left(1- \P\left(\max_{\t\in\I_{\T^{a}}} -\eta(\t)\leq -v_{\T}^{'}(x_{1}), \max_{\t\in \I_{\T^{a}}\cap\prod_{i=1}^{d} \Re(\delta_{i})} \eta(\t)\leq  v_{\T}^{*}(y_{1}) \right) \right) \nonumber\\
&&\sim e^{x_{1}+r+\sqrt{2r}z}+e^{y_{1}+r+\sqrt{2r}z}-H_{\a,\boldsymbol{\alpha}}^{-x_{1}+\log( \prod_{i=1}^{d}H_{\alpha_{i}}),-y_{1}+\log(\prod_{i=1}^{d} H_{a_{i},\alpha_{i}})}e^{r+\sqrt{2r}z}.
\end{eqnarray}
Following the same arguments as used in \eqref{eq8.2}, we have
\begin{eqnarray}
\label{eq9.5}
&&\prod_{i=1}^{d} n_{i} \left(\P\left(\max_{\t\in\I_{\T^{a}}} -\eta(\t)> -v_{\T}^{'}(x_{1}), \max_{\t\in \I_{\T^{a}}\cap\prod_{i=1}^{d} \Re(\delta_{i})} \eta(\t)>  u_{\T}^{*}(y_{2}) \right)\right) \to 0
\end{eqnarray}
as $ \T\to\infty $, and
\begin{eqnarray}
\label{eq9.6}
&&\prod_{i=1}^{d} n_{i} \left(\P\left(\max_{\t\in\I_{\T^{a}}} \eta(\t)> u_{\T}^{'}(x_{2}), \max_{\t\in \I_{\T^{a}}\cap\prod_{i=1}^{d} \Re(\delta_{i})} \eta(\t)\geq  -v_{\T}^{*}(y_{1}) \right)\right) \to 0
\end{eqnarray}
as $ \T\to\infty $. Combining \eqref{eq8.1}, \eqref{eq8.6}, \eqref{eq9.3}-\eqref{eq9.6} and Lemmas \ref{lem4}-\ref{lem7}, the assertion of this lemma follows.\qed

\section{Proofs of the main results}\label{sec4}
\noindent {\bf Proof of Theorem \ref{th1}.}
\pzx{Using Lemma \ref{lem8} and letting $ x_{1} \to -\infty$, $ y_{1} \to -\infty$, respectively,} we have
\begin{eqnarray}
\label{eq10.3}
&&\lim_{\T\to\infty}\P\left(\mathbf{M}_{\T}\leq \mathbf{u}_{\T},m_{\T}^{\boldsymbol{\delta}}>v_{\T}^{\boldsymbol{\delta}}(y_{1})\right)\nonumber\\
&=&\int_{-\infty}^{+\infty} \exp\Big( -(e^{-x_{2}-r+\sqrt{2r}z}
+e^{-y_{2}-r+\sqrt{2r}z}+e^{y_{1}+r+\sqrt{2r}z})\Big) \phi(z) dz,
\end{eqnarray}
and
\begin{eqnarray}
\label{eq10.4}
&&\lim_{\T\to\infty}\P\left(\mathbf{M}_{\T}\leq \mathbf{u}_{\T},m_{\T}>v_{\T}(x_{1})\right)\nonumber\\
&=&\int_{-\infty}^{+\infty} \exp\Big( -(e^{-x_{2}-r+\sqrt{2r}z}
+e^{-y_{2}-r+\sqrt{2r}z}+e^{x_{1}+r+\sqrt{2r}z})\Big) \phi(z) dz.
\end{eqnarray}

\pzx{Similarly, letting $ x_{1} \to -\infty$, $ y_{1} \to -\infty$ in the same time, we can get
\begin{eqnarray}
\label{eq10.41}
\lim_{\T\to\infty}\P\left(\mathbf M_{\T}\leq\mathbf u_{\T}\right)=\int_{-\infty}^{+\infty} \exp\Big( -(e^{-x_{2}-r+\sqrt{2r}z}
+e^{-y_{2}-r+\sqrt{2r}z})\Big) \phi(z) dz.
\end{eqnarray}
}

Therefore by \eqref{eq10.3}-\eqref{eq10.41} and Lemma \ref{lem8}, we can get the desired result  since
\begin{eqnarray}
\label{eq10.5}
\P\left(\mathbf M_{\T}\leq\mathbf u_{\T},\mathbf m_{\T}\leq\mathbf v_{\T}\right) &=& \P\left(\mathbf{v}_{\T}<\mathbf{m}_{\T}\leq \mathbf{M}_{\T}\leq \mathbf{u}_{\T}\right)-\P\left(\mathbf M_{\T}\leq\mathbf u_{\T}, m_{\T}> v_{\T}(x_{1})\right)\nonumber\\
&-&\P\left(\mathbf M_{\T}\leq\mathbf u_{\T}, m_{\T}^{\boldsymbol{\delta}}> v_{\T}^{\boldsymbol{\delta}}(y_{1})\right)+\P\left(\mathbf M_{\T}\leq\mathbf u_{\T}\right).
\end{eqnarray}
\qed

\noindent {\bf Proof of Theorem \ref{th2}.} It follows from Lemma \ref{lem9} and the definition of $ H_{\a,\boldsymbol{\alpha}}^{x,y} $ that
\begin{eqnarray*}
\lim_{\T\to\infty}\P\left(\mathbf{M}_{\T}\leq \mathbf{u}_{\T},m_{\T}>v_{\T}(x_{1})\right)&=&\int_{-\infty}^{+\infty} \exp\Big( -(e^{-x_{2}-r+\sqrt{2r}z}
+e^{-y_{2}-r+\sqrt{2r}z}\nonumber\\&-&H_{\a,\boldsymbol{\alpha}}^{x_{2}+\log( \prod_{i=1}^{d}H_{\alpha_{i}}),y_{2}+\log(\prod_{i=1}^{d} H_{a_{i},\alpha_{i}})}e^{-r+\sqrt{2r}z}+e^{x_{1}+r+\sqrt{2r}z})\Big) \phi(z) dz.
\end{eqnarray*}
Similarly,
\begin{eqnarray*}
\lim_{\T\to\infty}\P\left(\mathbf{M}_{\T}\leq \mathbf{u}_{\T},m_{\T}^{\boldsymbol{\delta}}>v_{\T}^{\boldsymbol{\delta}}(y_{1})\right)&=&\int_{-\infty}^{+\infty} \exp\Big( -(e^{-x_{2}-r+\sqrt{2r}z}
+e^{-y_{2}-r+\sqrt{2r}z}\nonumber\\&-&H_{\a,\boldsymbol{\alpha}}^{x_{2}+\log( \prod_{i=1}^{d}H_{\alpha_{i}}),y_{2}+\log(\prod_{i=1}^{d} H_{a_{i},\alpha_{i}})}e^{-r+\sqrt{2r}z}+e^{y_{1}+r+\sqrt{2r}z})\Big) \phi(z) dz,
\end{eqnarray*}
and
\pzx{
\begin{eqnarray*}
\lim_{\T\to\infty}\P\left(\mathbf{M}_{\T}\leq \mathbf{u}_{\T}\right)&=&\int_{-\infty}^{+\infty} \exp\Big( -(e^{-x_{2}-r+\sqrt{2r}z}
+e^{-y_{2}-r+\sqrt{2r}z}\nonumber\\&-&H_{\a,\boldsymbol{\alpha}}^{x_{2}+\log( \prod_{i=1}^{d}H_{\alpha_{i}}),y_{2}+\log(\prod_{i=1}^{d} H_{a_{i},\alpha_{i}})}e^{-r+\sqrt{2r}z})\Big) \phi(z) dz,
\end{eqnarray*}
}
By the same arguments as used in the proof of Theorem \ref{th1}, we can get the assertion of this theorem. \qed

\noindent {\bf Proof of Theorem \ref{th3}.} It follows from Lemma \ref{lem3} that
\begin{eqnarray*}
&&\Big|\prod_{i=1}^{d} n_{i} \left(1- \P\left(\max_{\t\in\I_{\T^{a}}} \eta(\t)\leq u_{\T}^{'}(x_{2}), \max_{\t\in \I_{\T^{a}}\cap\prod_{i=1}^{d} \Re(\delta_{i})} \eta(\t)\leq  u_{\T}^{*}(y_{2}) \right) \right) \nonumber\\
&&\qquad-\prod_{i=1}^{d} n_{i} \left(1- \P\left(\max_{\t\in\I_{\T^{a}}} \eta(\t)\leq u_{\T}^{'}(x_{2}), \max_{\t\in \I_{\T^{a}}} \eta(\t)\leq  u_{\T}^{*}(y_{2}) \right) \right)\Big| \nonumber\\
&\leq& \prod_{i=1}^{d} n_{i} \left(\P\left(\max_{\t\in \I_{\T^{a}}\cap\prod_{i=1}^{d} \Re(\delta_{i})} \eta(\t)\leq  u_{\T}^{*}(y_{2}) \right)-\P\left(\max_{\t\in\I_{\T^{a}}} \eta(\t)\leq u_{\T}^{*}(y_{2})\right)\right)
\to 0
\end{eqnarray*}
as $ \T\to\infty $. Then from \eqref{eq8.4} and \eqref{eq8.5}, we have
\pzx{
\begin{eqnarray}
\label{eq11.1}
&&\prod_{i=1}^{d} n_{i} \left(1-\P\left(\max_{\t\in\I_{\T^{a}}} \eta(\t)\leq u_{\T}^{'}(x_{2}), \max_{\t\in \I_{\T^{a}}\cap\prod_{i=1}^{d} \Re(\delta_{i})} \eta(\t)\leq  u_{\T}^{*}(y_{2}) \right)\right) \nonumber\\
&\to& e^{-\min\left(x_{2},y_{2}\right)-r+\sqrt{2r}z}
\end{eqnarray}
as $ \T\to\infty $, and
\begin{eqnarray}
\label{eq11.2}
&&\prod_{i=1}^{d} n_{i} \left(1-\P\left(\max_{\t\in\T_{\T^{a}}} -\eta(\t)<-v_{\T}^{'}(x_{1}),\max_{\t\in \I_{\T^{a}} \cap\prod_{i=1}^{d} \Re(\delta_{i})} -\eta(\t)< -v_{\T}^{*}(y_{1}) \right)\right)\nonumber\\
&\to& e^{\max\left(x_{1},y_{1}\right)+r+\sqrt{2r}z}
\end{eqnarray}
as $ \T\to\infty $.
Using Lemma \ref{lem3} again, we can get
\begin{eqnarray*}
&&\Big|\prod_{i=1}^{d} n_{i} \left( \P\left(\max_{\t\in\I_{\T^{a}}} \eta(\t)> u_{\T}^{'}(x_{2}), \max_{\t\in \I_{\T^{a}}\cap\prod_{i=1}^{d} \Re(\delta_{i})} -\eta(\t)\geq -v_{\T}^{*}(y_{1}) \right) \right) \nonumber\\
&&\qquad-\prod_{i=1}^{d} n_{i} \left( \P\left(\max_{\t\in\I_{\T^{a}}} \eta(\t)> u_{\T}^{'}(x_{2}), \max_{\t\in \I_{\T^{a}}} -\eta(\t)\geq  -v_{\T}^{*}(y_{1}) \right) \right)\Big| \nonumber\\
&\leq& \prod_{i=1}^{d} n_{i} \left(\P\left(\max_{\t\in\I_{\T^{a}}} -\eta(\t)\geq -v_{\T}^{*}(y_{1})\right)-\P\left(\max_{\t\in \I_{\T^{a}}\cap\prod_{i=1}^{d} \Re(\delta_{i})} -\eta(\t)\geq -v_{\T}^{*}(y_{1}) \right)\right)
\to 0
\end{eqnarray*}
as $ \T\to\infty $. Then by Borell theorem and \eqref{eq8.2}, we have
\begin{eqnarray}
\label{eq11.3}
&& \prod_{i=1}^{d} n_{i} \P\left(\max_{\t\in\I_{\T^{a}}} \eta(\t)> u_{\T}^{'}(x_{2}),\max_{\t\in \I_{\T^{a}} \cap\prod_{i=1}^{d} \Re(\delta_{i})} -\eta(\t)\geq -v_{\T}^{*}(y_{1}) \right)\to 0
\end{eqnarray}
as $ \T\to\infty $, and similarly,
\begin{eqnarray}
\label{eq11.4}
\prod_{i=1}^{d} n_{i} \P\left(\max_{\t\in\T_{\T^{a}}} -\eta(\t)\geq-v_{\T}^{'}(x_{1}), \max_{\t\in \I_{\T^{a}}\cap\prod_{i=1}^{d} \Re(\delta_{i})} \eta(\t)>  u_{\T}^{*}(y_{2})\right) \to 0
\end{eqnarray}
as $ \T\to\infty $.}

Hence, combining \eqref{eq8.1}, \eqref{eq8.6} and \eqref{eq11.1}-\eqref{eq11.4}, we have
\begin{eqnarray*}
\lim_{\T\to\infty}\P\left(\mathbf{v}_{\T}<\mathbf{m}_{\T}\leq \mathbf{M}_{\T}\leq \mathbf{u}_{\T}\right)=
\int_{-\infty}^{+\infty}\exp\Big(-e^{-\min(x_{2},y_{2})-r+\sqrt{2r}z}-e^{\max(x_{1},y_{1})+r+\sqrt{2r}z}\Big)\phi(z)dz.
\end{eqnarray*}
\pzx{By the same arguments as used in \eqref{eq10.3}-\eqref{eq10.41}}, we have
\begin{eqnarray*}
\lim_{\T\to\infty}\P\left(\mathbf{M}_{\T}\leq \mathbf{u}_{\T},m_{\T}>v_{\T}(x_{1})\right)=\int_{-\infty}^{+\infty} \exp\Big( -e^{-\min(x_{2},y_{2})-r+\sqrt{2r}z}
-e^{x_{1}+r+\sqrt{2r}z}\Big) \phi(z) dz,
\end{eqnarray*}

\begin{eqnarray*}
\lim_{\T\to\infty}\P\left(\mathbf{M}_{\T}\leq \mathbf{u}_{\T},m_{\T}^{\boldsymbol{\delta}}>v_{\T}(y_{1})\right)=\int_{-\infty}^{+\infty} \exp\Big( -e^{-\min(x_{2},y_{2})-r+\sqrt{2r}z}
-e^{y_{1}+r+\sqrt{2r}z}\Big) \phi(z) dz.
\end{eqnarray*}
and
\begin{eqnarray*}
\lim_{\T\to\infty}\P\left(\mathbf{M}_{\T}\leq \mathbf{u}_{\T}\right)=\int_{-\infty}^{+\infty} \exp\Big( -e^{-\min(x_{2},y_{2})-r+\sqrt{2r}z}\Big) \phi(z) dz,
\end{eqnarray*}
Therefore, from \eqref{eq10.5}, we complete the proof.
\qed


\begin{thebibliography}{999}

\bibitem{}
Adler, R. J. (2000). On excursion sets, tube formulas and maxima of random fields. {\it The Annals of Applied Probability}, {\bf 10(1)}, 1-74.

\bibitem{}
Adler, R. J. , Moldavskaya, E. , and Samorodnitsky, G. (2014). On the existence of paths between points in high level excursion sets of gaussian random fields. {\it The Annals of Probability}, {\bf 42(3)}, 1020-1053.


\bibitem{}
Berman, S. M. (1971). Asymptotic independence of the numbers of high and low level crossings of stationary gaussian processes. {\it Annals of Mathematical Statistics}, {\bf 42(3)}, 927-945.


\bibitem{}
Berman, S. M. (1974). Sojourns and extremes of gaussian processes. {\it Annals of Probability}, {\bf 2(6)}, 999-1026.

\bibitem{}
Chen, Y. and Tan, Z. (2016). Maxima and sum for discrete and continuous time Gaussian processes.
{\it Front Math China}, {\bf 11(1)}, 27-46.

\bibitem{}
Davis, R. A. (1979). Maxima and minima of stationary sequences. {\it Annals of Probability}, {\bf 7(3)}, 453-460.

\bibitem{}
D\c{e}bicki, K., Hashorva, E. and Soja-Kukiela. (2013). Extremes of homogeneous gaussian random fields. {\it Journal of Applied Probability}, {\bf 52(1)}, 55-67.

\bibitem{}
Hashorva, E. and Ji, L. (2016). Extremes of $ \alpha(t)$-locally stationary Gaussian random fields. {\it Transactions of the American Mathematical Society}, {\bf 368(1)}, 1-26.


\bibitem{}
Hashorva, E. and Tan, Z. (2015). Piterbarg's max-discretization theorem for stationary vector gaussian processes observed on different grids. {\it Statistics}, {\bf 49(2)}, 338-360.


\bibitem{}
H\"{u}sler, J. (2004). Dependence between extreme values of discrete and continous time locally stationary Gaussian processes.
{\it Extremes}, {\bf 7}, 179-190.

\bibitem{}
Leadbetter, M. R., Lindgren, G. and Rootz\'{e}n, H. (1983). {\it
Extremes and related properties of random sequences and processes}.
 Springer Verlag, New York.

\bibitem{}
Li, W. and Shao, Q. (2002). A normal comparison inequality and its applications. {\it Probability Theory and Related Fields}
{\bf 122(4)}, 494¨C508.


\bibitem{}
Liao, X. and Peng, Z. (2015). Asymptotics for the maxima and minima of h¨¹sler-reiss bivariate gaussian arrays. {\it Extremes}, {\bf 18(1)}, 1-14.

\bibitem{}
Lu, Y. and Peng, Z. (2017). Maxima and minima of independent and non-identically distributed
bivariate Gaussian triangular arrays. {\it Extremes}, {\bf 20(1)}, 187-198

\bibitem{}
Mccormick, W. P. and Qi, Y. (2000). Asymptotic distribution for the sum and maximum of gaussian processes. {\it Journal of Applied Probability}, {\bf 37(4)}, 958-971.


\bibitem{}
Mittal, Y. and Ylvisaker, D. (1975). Limit distributions for the maxima of stationary gaussian processes. {\it Stochastic Processes and Their Applications}, {\bf 3(1)}, 1-18.


\bibitem{}
Pickands, J. III. (1969). Asymptotic properties of the maximum in a stationary Gaussian procrss.
{\it Trans Amer Math Soc}, {\bf 145}, 75-86.


\bibitem{}
Piterbarg, V.I. (1996). {\it Asymptotic methods in the theory of Gaussian processes and fieldes}. AMS, Providence.

\bibitem{}
Piterbarg, V.I. (2004). Discrete and continuous time extremes of Gaussian processes.
{\it Extremes}, {\bf 7}, 161-177.

\bibitem{}
Peng, Z., Weng, Z., and Nadarajah, S. (2011). Joint limiting distributions of maxima and minima for complete and incomplete samples from weakly dependent stationary sequences. {\it Journal of Computational Analysis and Applications}, {\bf 13(5)}, 875-880.


\bibitem{}
\pzx{Song, Y., Turner, J. A., Peng, Z., Chen, C. and Li, X. (2018). Enhanced ultrasonic flaw detection using an
ultrahigh gain and time-dependent threshold. {\it IEEE Transactions on Ultrasonics, Ferroelectrics, and Frequency Control}, {\bf 65(7)}, 1214--1225.}

\bibitem{}
\pzx{Song, Y.,   Kube, C. M., Peng, Z., Turner, J. A. and Li, X. (2019). Flaw detection with ultrasonic backscatter signal
envelopes. {\it Journal of the Acoustical Society of America},  {\bf 145(2)},  142--148.
}

\bibitem{}
Turkman, K.F. (2012). Discrete and continuous time extremes of stationary processes.
{\it In: Rao, T., Rao, S. Rao, C., eds. Handbook of Statistics, Vol 30. Time Series Methods and Applications. Amsterdam: Elsevier}, 565-580.

\bibitem{}
Tan, Z. and Hashorva, E. (2014). \pzx{On Piterbarg's max-discretisation theorem for multivariate stationary Gaussian processes}. {\it Journal of Mathematical Analysis and Applications}, {\bf 409(1)}, 299-314.


\bibitem{}
Tan, Z. and Tang, L. (2014). The dependence of extremes values of discrete and continuous time strongly dependent Gaussian processes.
{\it Stochastics}, {\bf 86}, 60-69.


\bibitem{}
Tan, Z. and Wang, K. (2015). On Piterbarg's max-discretisation theorem for homogeneous Gaussian random fields. {\it Journal of Mathematical Analysis and Applications}, {\bf 429(2)}, 969-994.


\bibitem{}
Tan, Z. and Wang, Y. (2013). Extremes values of discrete and continuous time strongly dependent Gaussian processes.
{\it Communications in Statistics--Theory and Methods}, {\bf 42}, 2451-2463.








\end{thebibliography}
\end{document}